\documentclass[12pt]{article}

\usepackage{latexsym}
\usepackage{graphicx}
\usepackage{color}
\usepackage{multirow}

\usepackage{amsmath}
\usepackage{amsfonts}
\usepackage{amssymb}
\usepackage{graphicx}
\usepackage{float}
\usepackage{mathrsfs}

\graphicspath{{maxwell_images_png/}}
\usepackage[hidelinks]{hyperref}

\begin{document}

\title{Processing the 2D and 3D Fresnel experimental databases via topological derivative methods}

\author{A. Carpio, Universidad Complutense de Madrid, Spain, \\  
M. Pena, Universidad Polit\'ecnica de Madrid,  Spain, \\  
 M.L. Rap\'un, Universidad Polit\'ecnica de Madrid,  Spain}
 
 \date{May 24, 2021}

\maketitle

{\bf Abstract.}
This paper presents reconstructions of homogeneous targets from the 2D and 3D
Fresnel databases by one-step imaging methods based on the computation of topological derivative and topological energy fields. The electromagnetic inverse
scattering problem is recast as a constrained optimization problem, in which we seek
to minimize the error when comparing experimental microwave measurements with  computer-generated synthetic data for arbitrary targets by approximating
a Maxwell forward model. The true targets are then characterized by combining the topological derivatives or energies of such shape functionals for all available receivers
and emitters at different frequencies. Our approximations are comparable to the best approximations already obtained by other methods. However, these topological fields admit easy to evaluate closed-form expressions, which speeds up the process. \\

{\bf Keywords:} Microwave imaging, topological derivative, topological energy, multifrequency, Institut Fresnel databases, experimental data.

\section{Introduction}
\label{sec:intro}

Open access experimental datasets give researchers the opportunity of testing
inversion schemes against reliable data. Such is the case of the 2D and 3D datasets collecting microwave measurements in specific experimental setups developed
by the Institut Fresnel in Marseille, France. This initiative fostered several special
sections of {\it Inverse Problems} on ‘Testing inversion algorithms against
experimental data’, containing a series of contributions in which a variety of
inversion schemes were evaluated \cite{belkebir2001testing, geffrin2009continuing}.
The usefulness of such databases is twofold. On the one hand, they allow research
groups lacking experimental installations to test their methods against experimental
data. On the other  hand, they set a benchmark, such that different techniques can
be compared.

The setups considered in the Fresnel datasets launch microwaves towards targets made either of dielectric or conducting materials. The resulting electric field is
measured at a network of antennas. Our goal here is to test on these datasets topological field based schemes, which have been recently employed with
success in a variety of related inverse scattering problems
\cite{carpio2012hybrid,  funes2016defect,
masmoudi2005maxwell, Novotny, park1},
though often tested on synthetic data  (namely, numerically
generated). Few papers process actual experimental data via topological
derivative methods, see  \cite{dominguez2010non} for ultrasonic
applications,  \cite{tokmashev2013experimental} for elastic imaging,
\cite{yuan2013application} for vibroacoustic data, \cite{carpio2018experimental}
for holography and \cite{xavier2017topological} for fracture mechanics
applications, for instance.

The first Fresnel database contains data obtained in a geometrical configuration that
can be modelled by a scalar two-dimensional (2D) Helmholtz problem \cite{belkebir2001testing}.  A number of approaches were tested on it with variable
results: Bayesian approaches \cite{baussard2001bayesian}, modified gradient schemes \cite{belkebir2001modified2,duchene2001inversion,marklein2001linear},
the Born method \cite{belkebir2001modified2}, contrast source inversion techniques \cite{bloemenkamp2001inversion,marklein2001linear},
nonlinear inversion schemes  \cite{crocco2001inverse},
image fusion approaches \cite{fatone2001image},
linear diffraction tomography and real-coded genetic algorithms
\cite{marklein2001linear},
linear spectral estimation techniques \cite{testorf2001imaging},
multiple-frequency distorted-wave Born approaches \cite{tijhuis2001multiple},
and level-set based shape identification methods \cite{ramananjaona2001shape}.

In the three-dimensional (3D) database \cite{geffrin2009continuing} neither the objects
nor the geometrical arrangement display symmetries that allow for a simplification, therefore the full 3D Maxwell model has to be employed. This problem
is much more demanding as it can be concluded from  the variability of the results obtained by the algorithms initially tested against it: preliminary support reconstruction \cite{catapano20093d}, conjugate gradient-coupled dipole method \cite{chaumet2009three}, multiplicative smoothing with value picking regularization \cite{de2009three},  Bayesian framework with realistic random noise \cite{eyraud2009microwave}, multiplicative regularized contrast source inversion
method \cite{li2009application},  and a DBIM-BCGS method \cite{yu2009reconstruction}.


We will see that topological methods provide a useful non-iterative technique to
study these datasets. The paper is organized as follows. Section \ref{sec:framework}
recalls the mathematical model for the inverse scattering problem, while Section \ref{sec:td} describes the topological imaging tools employed.
Section \ref{sec:2D} presents reconstructions of objects obtained by topological
derivatives and topological energies from the 2D datasets. The results are comparable
to the best results obtained in previous studies by other methods, while the
computational complexity is very low. Section \ref{sec:3D} discusses the results
we have obtained for the 3D datasets.  These two sections follow a similar structure.
First we describe the geometrical configuration of the objects and the emitting and receiving antennas. Next, we adapt topological imaging methods to the setup under study and present the numerical results.
Finally, Section \ref{sec:conclusions} states our conclusions.


\section{Mathematical framework}
\label{sec:framework}

As explained in \cite{belkebir2001testing, geffrin2009continuing} the experiments
are carried out in an anechoic chamber, hence the physical model can be posed in
${\cal R}=\mathbb{R}^3$. We denote by $\Omega\subset{\cal R}$ the region occupied
by the objects and assume that the ambient medium surrounding the objects has the same properties as vacuum. When the incident wave field is time harmonic, that is, it takes the form ${\bf E}^{\rm inc}({\bf x},t)= \mbox{Re}\left(e^{-\imath \omega t} \boldsymbol{\mathcal{E}}^{\rm inc}(\mathbf x)\right)$,
the wave field solution of Maxwell's equations will be time harmonic too.
The complex amplitude
$\boldsymbol{\mathcal{E}}(\mathbf x)$ of time harmonic electromagnetic waves
${\bf E}({\bf x},t)=\mbox{Re}\left(e^{-\imath \omega t}
 \boldsymbol{\mathcal{E}}(\mathbf x)\right)$ satisfies
\[
\mathbf{curl}\, \mathbf{curl}\left(\boldsymbol{\mathcal{E}}\right)- \kappa^2\boldsymbol{\mathcal{E}}=\mathbf{0}\quad\mathrm{in}\,{\cal R}\setminus\overline{\Omega}.
\]
Here, $\kappa =  \omega\sqrt{\mu_0 \varepsilon_0}$ is the wave number,
$\mu_0 $ and $\varepsilon_0$ being the permeability and permittivity of vacuum,
while $\omega$  represents the angular frequency,  related to the
frequency $\nu$ by $\omega = 2 \pi \nu$.  This equation is complemented by the
Silver-Müller condition at infinity \cite{monk2003finite}
\begin{equation}\label{eq:SilverMuller}
\lim_{|\mathbf x| \rightarrow \infty} |\mathbf x|  \left|
\mathbf{curl} \,  (\boldsymbol{\mathcal{E}} - \boldsymbol{\mathcal{E}}^{\rm inc}) \times {\mathbf x \over |\mathbf x|}
-\imath \kappa (\boldsymbol{\mathcal{E}}- \boldsymbol{\mathcal{E}}^{\rm inc}) \right| =0.
\end{equation}
The corresponding equation for the magnetic field is omitted as only the electric
field was measured.

Two types of homogeneous targets are considered \cite{belkebir2001testing, geffrin2009continuing}. On  the one hand, dielectric
objects $\Omega$ with relative permeability $\mu_\mathrm{r}=1$ and different
values for the relative permittivity $\varepsilon_\mathrm{r}=\varepsilon_\mathrm{d}$.
In this case, the  conductivity $\sigma=0$ everywhere and the complex amplitude
$\boldsymbol{\mathcal{E}}$ satisfies
\[
\mathbf{curl}\, \mathbf{curl}\left(\boldsymbol{\mathcal{E}}\right)-\varepsilon_\mathrm{d}\kappa^2\boldsymbol{\mathcal{E}}=\mathbf{0}\quad\mathrm{in}\; \Omega,
\]
with transmission conditions at the interface
\[
\boldsymbol{\mathcal{E}}^+\times\mathbf{n}-\boldsymbol{\mathcal{E}}^-\times\mathbf{n}=\mathbf{0}, \quad
\mathbf{curl}\left(\boldsymbol{\mathcal{E}}\right)^+\times\mathbf{n}-\mathbf{curl}\left(\boldsymbol{\mathcal{E}}\right)^-\times\mathbf{n}=\mathbf{0}\quad\mathrm{on}\;\partial\Omega,
\]
where we have used the fact that $\mu=\mu_0$ everywhere. Here, $\mathbf n$
represents the  outer normal vector on the object's surface. The
symbols $+$ and $-$ denote limits from the exterior and interior,  respectively.

On the other hand,  we consider perfectly conducting targets for which
$\sigma\to\infty$. Then, Ohm's law  implies
$\boldsymbol{\mathcal{E}}=\mathbf{0} \; \mathrm{in}\; \Omega $
and the boundary condition at the interface becomes
$\boldsymbol{\mathcal{E}}\times\mathbf{n}=\mathbf{0}\;
\mathrm{on}\; \partial\Omega. $

Summarizing, for the three dimensional case, the electric field  will be modelled by
the transmission problem

\begin{equation}\label{eq:3ddiel}
\left\{
\begin{array}{ll}
\mathbf{curl}\; \mathbf{curl}\left(\boldsymbol{\mathcal{E}}\right)-\kappa^2\boldsymbol{\mathcal{E}}=\mathbf{0}, &\quad\mathrm{in}\,{\cal R}\setminus\overline{\Omega},\\
\mathbf{curl}\; \mathbf{curl}\left(\boldsymbol{\mathcal{E}}\right)-\varepsilon_\mathrm{d}\kappa^2\boldsymbol{\mathcal{E}}=\mathbf{0}, &\quad\mathrm{in}\,\Omega,\\
\boldsymbol{\mathcal{E}}^+\times\mathbf{n}-\boldsymbol{\mathcal{E}}^-\times\mathbf{n}=\mathbf{0}, &\quad\mathrm{on}\;\partial\Omega,\\
\mathbf{curl}\left(\boldsymbol{\mathcal{E}}\right)^+\times\mathbf{n}-\mathbf{curl}\left(\boldsymbol{\mathcal{E}}\right)^-\times\mathbf{n} =\mathbf{0},&\quad\mathrm{on}\;\partial\Omega,\\
\lim_{|\mathbf x| \rightarrow \infty} |\mathbf x|  \left|
\mathbf{curl} \,  (\boldsymbol{\mathcal{E}} - \boldsymbol{\mathcal{E}}^{\rm inc}) \times {\mathbf x \over |\mathbf x|}
-\imath \kappa (\boldsymbol{\mathcal{E}}- \boldsymbol{\mathcal{E}}^{\rm inc}) \right| =0,
\end{array}
\right.
\end{equation}
in the experiments with dielectric targets, whereas it will take the form of an exterior Dirichlet  problem
\begin{equation}\label{eq:3dcond}
\left\{
\begin{array}{ll}
\mathbf{curl}\; \mathbf{curl}\left(\boldsymbol{\mathcal{E}}\right)-\kappa^2\boldsymbol{\mathcal{E}}=\mathbf{0}, &\quad\mathrm{in}\,{\cal R}\setminus\overline{\Omega},\\
\boldsymbol{\mathcal{E}}\times\mathbf{n}=\mathbf{0}, &\quad\mathrm{on}\;\partial\Omega,\\
\lim_{|\mathbf x| \rightarrow \infty} |\mathbf x|  \left|
\mathbf{curl} \,  (\boldsymbol{\mathcal{E}} - \boldsymbol{\mathcal{E}}^{\rm inc}) \times {\mathbf x \over |\mathbf x|}
-\imath \kappa (\boldsymbol{\mathcal{E}}- \boldsymbol{\mathcal{E}}^{\rm inc}) \right| =0,
\end{array}
\right.
\end{equation}
for the conducting targets.

In the 2D Fresnel dataset \cite{belkebir2001testing}, the targets $\Omega$ are vertical cylinders, that is,
$\Omega=\Omega_\mathrm{2D}\times\left[-L_z,L_z\right], $
where $\Omega_\mathrm{2D}\subset\mathbb{R}^2$ represents the projection of the targets on the horizontal plane and \mbox{$L_z\gg \sqrt{\mathrm{Area}\left(\Omega_\mathrm{2D}\right)}$}. Since the  length along the axis direction is much larger than the horizontal characteristic length, we can make the approximation
$\Omega\approx\Omega_\mathrm{2D}\times\mathbb{R}$ so that
$\frac{\partial \varepsilon_\mathrm{r}}{\partial z}=0,$  i.e., the relative electrical permittivity does not depend on the vertical coordinate. In this case, the Silver-Müller radiation condition at infinity (\ref{eq:SilverMuller}) is replaced by the Sommerfeld ratiation condition:
\begin{equation}\label{eq:Sommerfeld}
\lim_{\left|\mathbf{x}\right|\to\infty}\sqrt{\left|\mathbf{x}\right|}\left(\frac{\partial\left(\mathcal{E}-\mathcal{E}^\mathrm{inc}\right)}{\partial \left|\mathbf{x}\right|}-\imath\kappa\left(\mathcal{E}-\mathcal{E}^\mathrm{inc}\right)\right) =0.
\end{equation}
This condition ensures  that the scattered field $\mathcal{E}-\mathcal{E}^\mathrm{inc}$ tends to zero fast enough as $ |{\bf x}| $ tends to infinity, and selects only outward radiating solutions. Let us stress that this condition depends on the structure of the time harmonic waves. If we had chosen $ {\bf E}^{\rm inc}({\bf x},t)=\rm{Re}\left(e^{\imath \omega t} \boldsymbol{\mathcal E}^\mathrm{inc}(\mathbf x)\right)$, the sign of $\imath \kappa$  in the Sommerfeld condition would be positive. In fact, that is the convention used in the 2D Fresnel database, so, we must conjugate their data to use it with our convention, which we chose for consistency with the previous 3D case.

If we write the $\mathbf{curl} \, \mathbf{curl}$ operator in terms of laplacians and
divergences, we find  that for dielectric targets, the first two equations in
(\ref{eq:3ddiel}) can be written as:
\begin{equation*}
\left\{
\begin{array}{ll}
\mathrm{div}\left(\boldsymbol{\mathcal{E}}\right)=0, \quad
\Delta\boldsymbol{\mathcal{E}}+\kappa^2\boldsymbol{\mathcal{E}}=\mathbf{0},&\quad\mathrm{in}\,{\cal R}\setminus\overline{\Omega},\\
\mathrm{div}\left(\boldsymbol{\mathcal{E}}\right)=0, \quad
\Delta\boldsymbol{\mathcal{E}}+\varepsilon_\mathrm{d}\kappa^2\boldsymbol{\mathcal{E}}=\mathbf{0},&\quad\mathrm{in}\,\Omega.
\end{array}
\right.
\end{equation*}
This system allows for solutions of the form
$\boldsymbol{\mathcal{E}}\left(x,y,z\right)=\mathcal{E}\left(x,y\right)\mathbf{k},$
and problem (\ref{eq:3ddiel}) reduces to
\begin{equation} \label{eq:2ddiel}
\left\{
\begin{array}{ll}
\Delta_\mathrm{2D}\mathcal{E}+\kappa^2\mathcal{E}=0,&\quad\mathrm{in}\,{\cal R}_\mathrm{2D}\setminus\overline{\Omega}_\mathrm{2D}, \\
\Delta_\mathrm{2D}\mathcal{E}+\varepsilon_\mathrm{d}\kappa^2\mathcal{E}=0,&\quad\mathrm{in}\,\Omega_\mathrm{2D},\\
\mathcal{E}^+-\mathcal{E}^-=0,&\quad\mathrm{on}\;\partial\Omega_\mathrm{2D}, \\
\nabla_\mathrm{2D}\mathcal{E}^+\cdot\mathbf{n}_\mathrm{2D}-\nabla_\mathrm{2D}\mathcal{E}^-\cdot\mathbf{n}_\mathrm{2D}=0,&\quad\mathrm{on}\;\partial\Omega_\mathrm{2D},\\
\lim_{\left|\mathbf{x}\right|\to\infty}\sqrt{\left|\mathbf{x}\right|}\left(\frac{\partial\left(\mathcal{E}-\mathcal{E}^\mathrm{inc}\right)}{\partial \left|\mathbf{x}\right|}-\imath\kappa\left(\mathcal{E}-\mathcal{E}^\mathrm{inc}\right)\right) =0,
\end{array}
\right.
\end{equation}
where $\Delta_\mathrm{2D}$ and $\nabla_\mathrm{2D}$ are the horizontal part of the laplacian and the gradient:
\[\Delta_\mathrm{2D}:f\mapsto\frac{\partial f}{\partial x^2}+\frac{\partial f}{\partial y^2},
\quad
\nabla_\mathrm{2D}:\mathbf{f}\mapsto\frac{\partial f}{\partial x}\mathbf{i}+\frac{\partial f}{\partial y}\mathbf{j},\]
$\mathbf{n}_\mathrm{2D}$ is the two dimensional  outer normal to $\Omega_\mathrm{2D}$ and ${\cal R}_\mathrm{2D}=\mathbb R^2.$
Similarly, we have a two dimensional reduction for the conducting targets:
\begin{equation} \label{eq:2dcond}
\left\{
\begin{array}{ll}
\Delta_\mathrm{2D}\mathcal{E}+\kappa^2\mathcal{E}=0,&\quad\mathrm{in}\,{\cal R}\setminus\overline{\Omega}_\mathrm{2D},\\
\mathcal{E}^+ =0,&\quad\mathrm{on}\;\partial\Omega_\mathrm{2D},\\
\lim_{\left|\mathbf{x}\right|\to\infty}\sqrt{\left|\mathbf{x}\right|}\left(\frac{\partial\left(\mathcal{E}-\mathcal{E}^\mathrm{inc}\right)}{\partial \left|\mathbf{x}\right|}-\imath\kappa\left(\mathcal{E}-\mathcal{E}^\mathrm{inc}\right)\right) =0.
\end{array}
\right.
\end{equation}
In both cases we obtain 2D scalar Helmholtz problems. For simplicity,
we will omit the 2D subindex when we discuss the two dimensional results.

Problems (\ref{eq:3ddiel}), (\ref{eq:3dcond}), (\ref{eq:2ddiel}) and (\ref{eq:2dcond}),
constitute the forward models governing the  wave fields interacting with the
targets. Evaluating their solutions at the selected receivers, we should obtain the
data recorded for the true objects, except for perturbations due to experimental noise. In the
next section, we describe the topological derivative based approach to recover
the original targets from the measurements.

\section{Topological sensitivity based imaging}
\label{sec:td}

The inverse scattering problem underlying the Fresnel databases consists in finding
objects $\Omega$ such that the solutions of the corresponding forward problems
(\ref{eq:3ddiel}), (\ref{eq:3dcond}), (\ref{eq:2ddiel}) or (\ref{eq:2dcond}) evaluated
at the receivers are close to the recorded data for the selected incident waves.
This kind of problems is severely ill-posed \cite{colton2011ipe} and are often regularized recasting them as constrained optimizations problems: Find $\Omega$ minimizing an error functional $\mathcal{J}(\Omega)$ which compares the recorded data with the synthetic data that would be obtained solving the forward problem with target $\Omega$.
A typical choice is
\begin{eqnarray}
\mathcal{J}(\Omega)  = {1\over 2} \sum_{j=1}^{N_{\rm meas}}
|\boldsymbol{\mathcal{E}}(\mathbf x_j) - \boldsymbol{\mathcal{E}}^{\rm meas}_j |^2,
\end{eqnarray}
where $\mathbf x_j$, $j=1,..., N_{\rm meas}$  are the receivers' positions
and $ \boldsymbol{\mathcal{E}}^{\rm meas}_j$  is the recorded data at the $jth-$receiver.
The topological derivative of such functionals has the potential of providing guesses
of the objects  \cite{feijoo2004td,masmoudi2005maxwell}, which can be very sharp
when enough frequencies and/or incident directions are combined
\cite{carpio2008topological,carpio2010determining,funes2016defect,guzina2015high,
park1,park2}.

Let $\mathcal{J}$ be a real valued shape functional and let
$\mathcal{R}_\epsilon=\mathcal{R}\setminus\overline{\mathrm{B}_\epsilon\left(\mathbf{x}\right)}$ be the domain obtained when a ball of radius $\epsilon$ centered at $\mathbf{x}$  is removed from  $\mathcal{R}$. The topological derivative is  defined \cite{sokolowski1999td} as a real scalar field $\mathrm{D}_\mathrm{T}\mathcal{J}:\mathcal{R} \to\mathbb{R}$ for which the following asymptotic expansion holds:
\begin{equation}\label{eq:shape_fun_ex}
\mathcal{J}\left(\mathcal{R}_\epsilon\right) = \mathcal{J} \left(\mathcal{R}\right) +\mathrm{D}_\mathrm{T}\mathcal{J}\left(\mathbf{x}\right) f(\epsilon)+
o\left(f(\epsilon)\right) \end{equation}
where $f(\epsilon)$ is a positive monotonically increasing function satisfying ${\lim_{\epsilon\to 0^+}f\left(\epsilon\right)=0}$ and selected to ensure that $\mathrm{D}_\mathrm{T}\mathcal{J}\left(\mathcal{R}\right)$ is finite and is not identically zero. Depending on the dimension of the problem and on the kind of targets, function $f$ should be selected as follows (see \cite{carpio2008topological, guzina2006misfit} for 2D dielectric objects, \cite{carpio2010determining, samet2003dirichlet} for 2D perfectly conducting ones, and \cite{le2017topological, masmoudi2005maxwell} for both 3D cases): for dielectric objects,
\[
f(\epsilon)=\left\{
\begin{array}{ll}
\pi\epsilon^2,&\quad \mbox{in 2D},\\
\frac43\,\pi\epsilon^3,&\quad\mbox{in 3D},
\end{array}
\right.
\]
while for perfectly conducting targets,
\[
f(\epsilon)=\left\{
\begin{array}{ll}
\frac{-2\pi}{\ln(\kappa\epsilon)},&\quad \mbox{in 2D},\\
\frac43\,\pi\epsilon^3,&\quad\mbox{in 3D}.
\end{array}
\right.
\]
Notice that for 2D perfectly conducting targets this function has not a clear geometrical interpretation, while in the remaining cases $f$ is the measure of the ball.

The topological derivative measures the sensitivity of the misfit functional to locating infinitesimal objects at points $\mathbf{x}$. Expansion (\ref{eq:shape_fun_ex}) indicates the cost functional decreases by locating objects at that points where the topological derivative is negative. Following that idea,  we will determine a family of domains at which the topological derivative takes large negative values:
\begin{equation}\label{eq:dom_app}
\Omega_{\mathrm{app},\lambda}=
\left\lbrace\mathbf{x}\in  {\cal R}^{\rm insp},\;\mathrm{s.t.}\;
\mathrm{D}_\mathrm{T}\mathcal{J}\left(\mathbf{x}\right)\le
\lambda\min_{\mathbf{y} \in {\cal R}^{\rm insp}}
\mathrm{D}_\mathrm{T}\mathcal{J}\left(\mathbf{y} \right)
\right\rbrace, \quad \lambda \in \left[0,1\right],
\end{equation}
which we expect to provide reconstructions of the true targets, where ${\cal R}^{\rm insp}$ is a big enough bounded subset of ${\cal R}$ where objects are assumed to be located.
We will consider our approximation method robust if the domains  $\Omega_{\mathrm{app},\lambda}$ remain close (with respect to some metric)
to the true scatterers $\Omega$ for a wide range of $\lambda$.

In practice, one can calculate  closed-form expressions for topological derivatives of many shape functionals in terms of solutions of forward and adjoint problems that take place in the pristine media. The peaks of a companion field, the topological energy, defined as the product of norms of forward and adjoint fields
\cite{dominguez2010non}, also mark the location of targets.
We will discuss these methods in more detail when applying them to the 2D and 3D
databases in the next two sections.

\section{Application to the two dimensional Fresnel database}
\label{sec:2D}

As already explained, in the 2D Fresnel database \cite{belkebir2001testing}, all objects are cylinders perpendicular to the plane where both the emitting and receiving antennas are placed (the horizontal plane). Antennas launch waves at a certain frequency
towards the target. Both the incident field, i.e. the total electric field when there is no object, as well as the total electric field when an object is present are measured at $N_\mathrm{meas}=49$ locations.
For each frequency,  $N_\mathrm{inc}=36$ experiments are performed by rotating the target from one experiment to the next up to $360^\circ$ degrees. In practice, we
assume that the object is stationary and the emitters/receivers rotate around it.
The emitting positions are
\begin{equation}\label{position}
\mathbf{x}^\mathrm{E} = R^\mathrm{E}\left( \cos\theta^\mathrm{E}\mathbf{i} + \sin\theta^\mathrm{E}\mathbf{j}\right)\quad\mathrm{with}\,\theta^\mathrm{E}=0^\circ,10^\circ,\dots,350^\circ,
\end{equation}
where $R^\mathrm{E}=0.76\,\mathrm{m}$ is the distance from the emitting antenna to the center of rotation of the target. The receiving positions are linearly spaced in the circumference  of radius $R^\mathrm{R}=0.72\,\mathrm{m}$ and separated from the emitting antenna $60^\circ$ degrees:
\[
\mathbf{x}^\mathrm{R} = R^\mathrm{R}\left( \cos\left(\theta^\mathrm{E}+\theta^\mathrm{R}\right)\mathbf{i}+ \sin\left(\theta^\mathrm{E}+\theta^\mathrm{R}\right)\mathbf{j}\right)\quad\mathrm{with}\,\theta^\mathrm{R}=60^\circ,65^\circ,\dots,300^\circ.
\]

\begin{figure}[h!]
\centering
\includegraphics[trim = 0cm 1cm 0cm 0.5cm, clip=true]{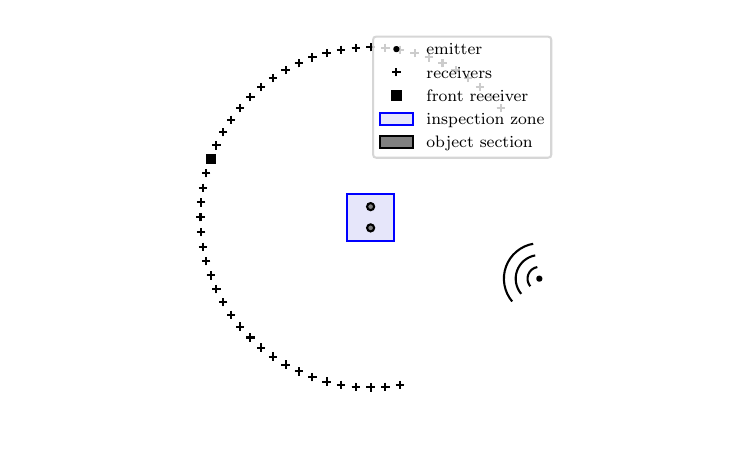}
\caption{Scaled representation of the location of the receiving and emitting antennas for a generic experiment, as well as the size of the inspection zone compared with the size of the  objects' section.}\label{fig:2d_diagram}
\end{figure}

\begin{figure}[h!]
\centering
\includegraphics[trim = 0cm 0cm 1cm 0cm, clip=true]{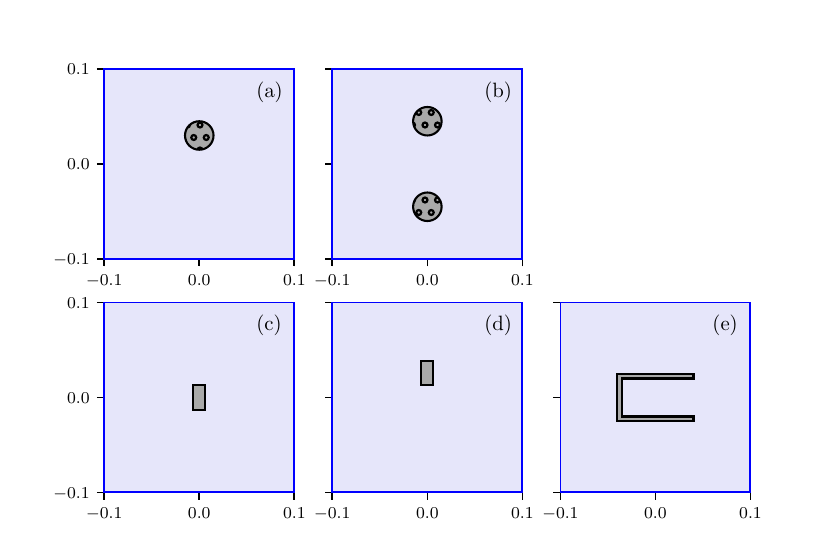}
\caption{Horizontal sections of the five different targets compared with the inspection zone. The first two ones are made of dielectric material (indicated by a dotted pattern) whereas the last three ones are metallic targets (homogeneous pattern).
}\label{fig:2dsections}
\end{figure}

Figure \ref{fig:2d_diagram} represents this setup for an arbitrary angle of the emitting antenna (angle of rotation of the target). The inspection zone  where we plot our reconstructions will be the square ${\cal R}^{\rm insp}=[-0.1,0.1]^2$. We will consider the data recorded in the transverse magnetic (TM) configuration, in which the polarization vector of the emitting and receiving antennas is vertical, i.e. $\mathbf{p}^E=\mathbf{p}^R=\mathbf{k}$ for all  the emissions $E$ and receptions $R$. This allows us to reduce the general 3D Maxwell model to a
scalar 2D Helmholtz model for the vertical component of the electric field.
Instead, in the transverse electric (TE) configuration the polarization vectors of the antennas lie in the horizontal plane. A 2D reduction would no longer be described by a scalar model, since the horizontal components of the electric field cannot be decoupled. Only one dataset corresponds to that configuration, and we will not study it here.
Figure \ref{fig:2dsections} depicts the 2D sections of the five considered targets.
Targets (a) and (b) are dielectric objects with $\varepsilon_\mathrm{d} = 3$, hence,
the total electric field is modelled by the transmission problem (\ref{eq:2ddiel}).
Targets (c), (d) and (e) are conducting objects, thus, the total electric field
is modelled by the exterior Dirichlet problem (\ref{eq:2dcond}). The ranges of frequencies
of the available datasets are collected in Table \ref{tab:exps}.

\begin{table}[h!]
\centering
\begin{tabular}{l|l|c|l|l}
Target & File name & Polarization &  Freq. band &  Freq. step \\
\hline
\hline
Single dielectric &  {\tt dielTM\_dec4f.exp} & TM & $4-16$  & $4$ \\
Single dielectric &  {\tt dielTM\_dec8f.exp} & TM & $1-8$  & $1$ \\
Two dielectrics &  {\tt twodielTM\_4f.exp} & TM & $4-16$  & $4$ \\
Two dielectrics &  {\tt twodielTM\_8f.exp} & TM & $1-8$  & $1$ \\
Metallic rectangle &  {\tt rectTM\_cent.exp} & TM & $4-16$  & $4$ \\
Metallic rectangle &  {\tt rectTM\_dece.exp} & TM & $2-16$  & $2$ \\
Metallic U &  {\tt uTM\_shaped.exp} & TM & $2-16$  & $2$
\end{tabular}
\caption{ Names of the original files in \cite{belkebir2001testing} and range of the frequencies  in GHz  used in each file. }\label{tab:exps}
\end{table}

The incident waves are not known explicitly, but through the values recorded
in the different experiments. Our first task will be to propose fittings for the
incident waves using the measured values. Then, we will obtain reconstructions
of the targets by topological methods.

\subsection{Fitting the incident waves}
\label{sec:incident2d}

\begin{figure}[h!]
\centering
\includegraphics[width=10cm]{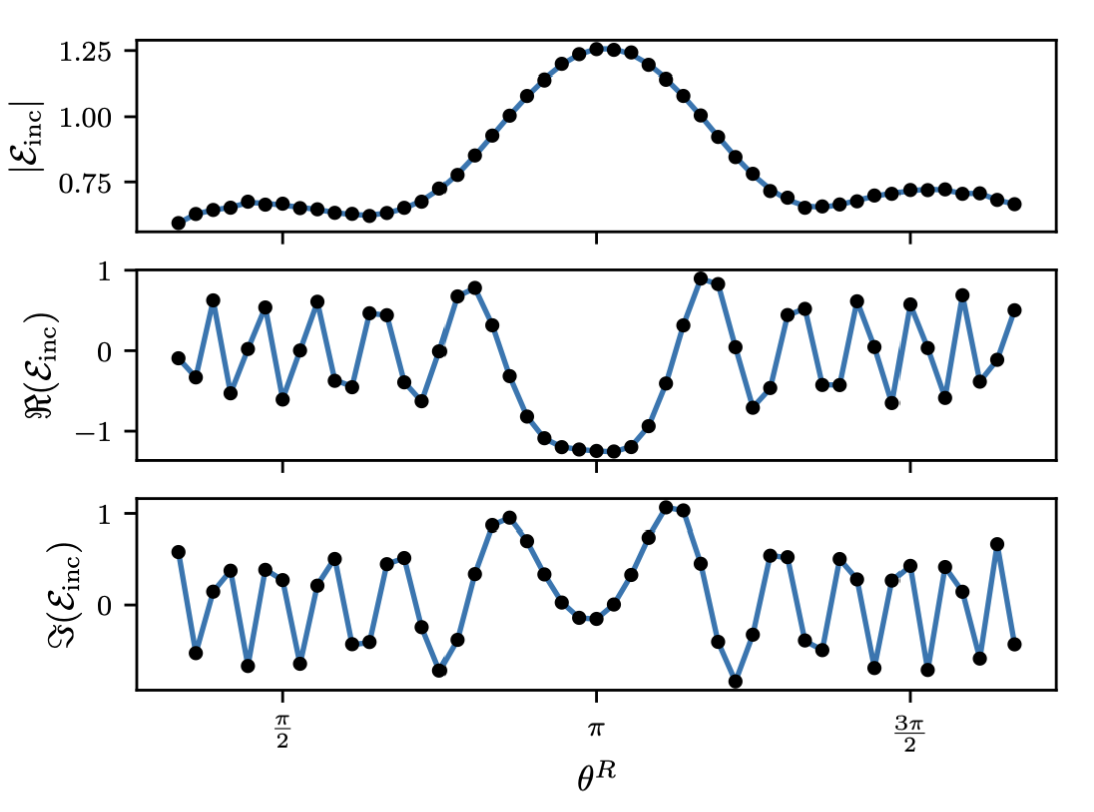}
\caption{Measured incident field (dots) for the experiments performed at a low frequency $2\,\mathrm{GHz}$ and fitted field (line). The amplitude is largest at the furthest receivers ($\theta^\mathrm{R}=\pi$) which suggests a highly directional antenna.}
\label{fig:3_fitted_low}
\end{figure}

Figure \ref{fig:3_fitted_low} shows the values of the incident field measured for $\nu = 2\,\mathrm{GHz}$ in the absence of targets. The highest amplitude is reached at $\theta^\mathrm{R}=\pi$, that is, at a point which is in front of the emitter. In order to model this directionality, we will fit a combination of  Hankel functions by the least-squares method. More precisely, we seek fittings of the form
\begin{eqnarray}
\mathcal{U}\left(\mathbf{x}\right)= a_0H^1_0\left( \kappa \rho\right)+ \hskip -2mm \sum_{n=1}^{N_\mathrm{modes}} \hskip -2mm \left(a_nH_n^1\left( \kappa \rho\right)\cos\left(n\theta\right)+b_nH_n^1\left( \kappa \rho\right)\sin\left(n\theta\right)\right)
\label{fitted}
\end{eqnarray}
where $H^1_n$ is the Hankel function of first kind and order $n$  (see \cite[Section 9]{abramowitz1948handbook}) and $\left(\rho,\theta\right)$ are polar coordinates
centered at the emitter:
\begin{eqnarray*}
\rho &=\left|\mathbf{x}-\mathbf{x}^\mathrm{E}\right|,\quad
\theta &= \arccos\left(\frac{\left(\mathbf{x}-\mathbf{x}^\mathrm{E}\right)\cdot
\mathbf{x}^\mathrm{E}}{|\mathbf{x}^\mathrm{E}|^2}\right).
\end{eqnarray*}

\begin{figure}[h!]
\centering
\includegraphics[width=10cm]{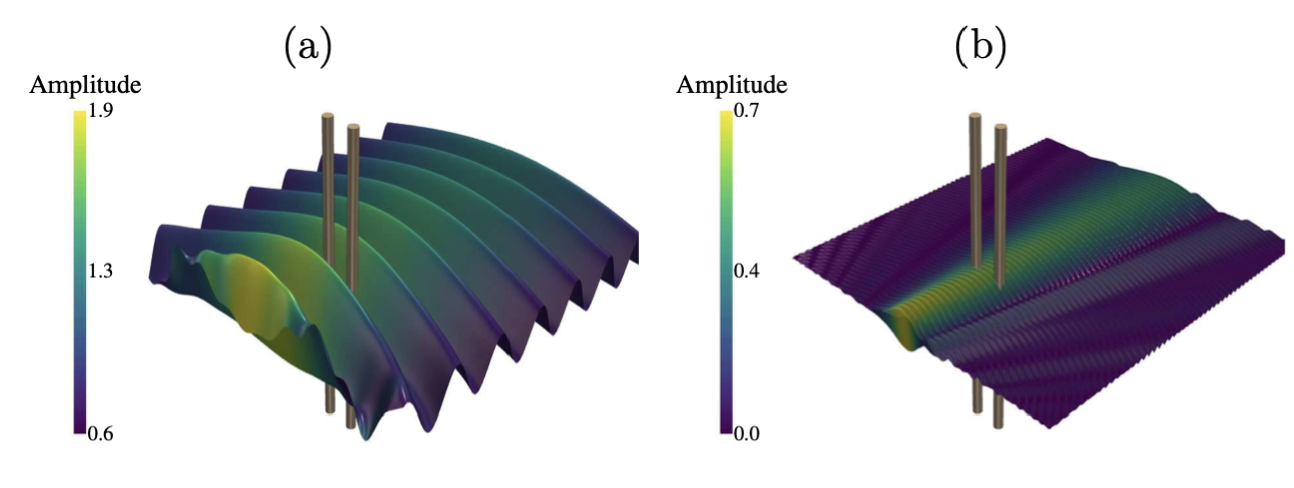}
\caption{(a) Incident field fitted by a series of Hankel functions for $2\,\mathrm{GHz}$. Notice the directionality of the antenna, as well as the resemblance of the incident wave
to a plane wave near the objects. (b) Incident field fitted  for the frequency  $16\,\mathrm{GHz}$. Now the field is far from isotropic.}
\label{fig:3d_2dwave_low}
\end{figure}

\noindent Figure \ref{fig:3_fitted_low} also represents the fitted incident wave for $N_\mathrm{modes}=14$. The fit is sharp. Plotting the real part of the fitted incident
field in the whole plane, we notice that anisotropy is less pronounced at the scale of
the targets, see Figure \ref{fig:3d_2dwave_low}(a).
This remark suggests  approximating the incident field as a plane wave or as an isotropic wave using only the field measured at the front receiver. The plane wave approximation would be
\begin{equation}\label{plane}
\mathcal{U}\left(\mathbf{x}\right)=\frac{\mathcal{E}^\mathrm{inc}_\mathrm{front}}{e^{i\left( \kappa \mathbf{d}\cdot\mathbf{x}^R_\mathrm{front}\right)}}e^{i \kappa \mathbf{d}\cdot\mathbf{x}}
\end{equation}
where $\mathbf{d}=\frac{\mathbf{x}_{\rm front}^R-\mathbf{x}^E}{|\mathbf{x}_{\rm front}^R-\mathbf{x}^E|}$ is the direction of propagation of the wave, and $\mathcal{E}^\mathrm{inc}_\mathrm{front}$
is the  incident field measured  in front of the emitter.  If we model the incident wave
as an isotropic wave, we avoid computing this direction vector  while obtaining similar values. In this case, we  approximate the incident wave  by
\begin{equation}\label{eq:state_planewave2d}
\mathcal{U}\left(\mathbf{x}\right)=\frac{\mathcal{E}^\mathrm{inc}_\mathrm{front}}{H_0^1\left( \kappa \left|\mathbf{x}^R_\mathrm{front}-\mathbf{x}^E\right|\right)}H_0^1\left( \kappa \left|\mathbf{x}-\mathbf{x}^E\right|\right).
\end{equation}

\begin{figure}[h!]
\centering
\includegraphics[width=10cm]{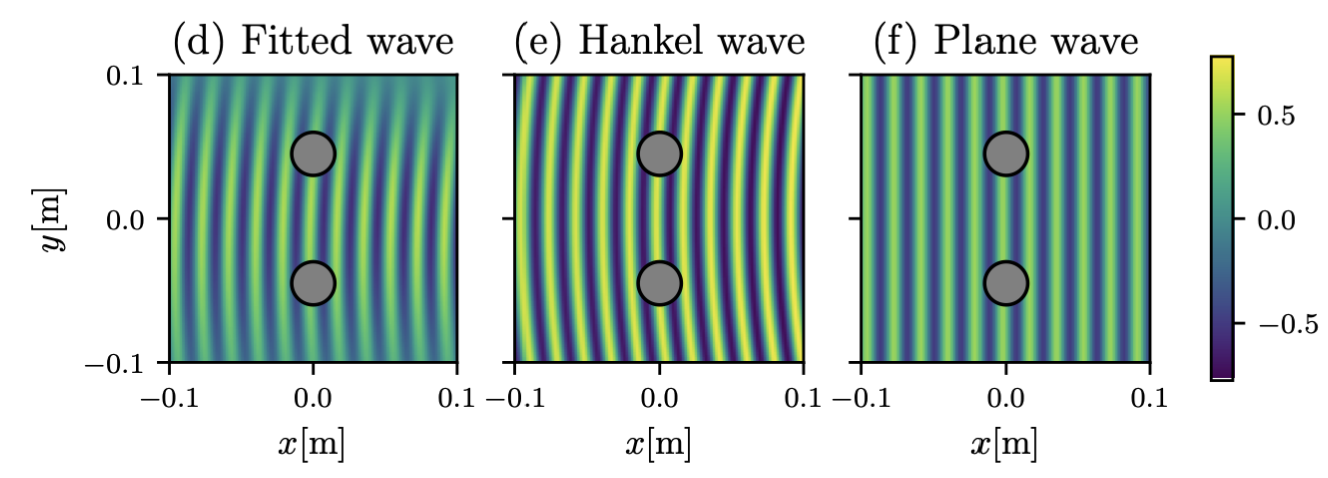}
\caption{(a)-(c) Incident electric field for the three wave models at
$2\,\mathrm{GHz}$.
(d)-(f) Incident electric field for the three wave models at $16\,\mathrm{GHz}$.
}
\label{fig:2d_comparisons_low}
\end{figure}

\begin{figure}[h!]
\centering
\includegraphics[width=10cm]{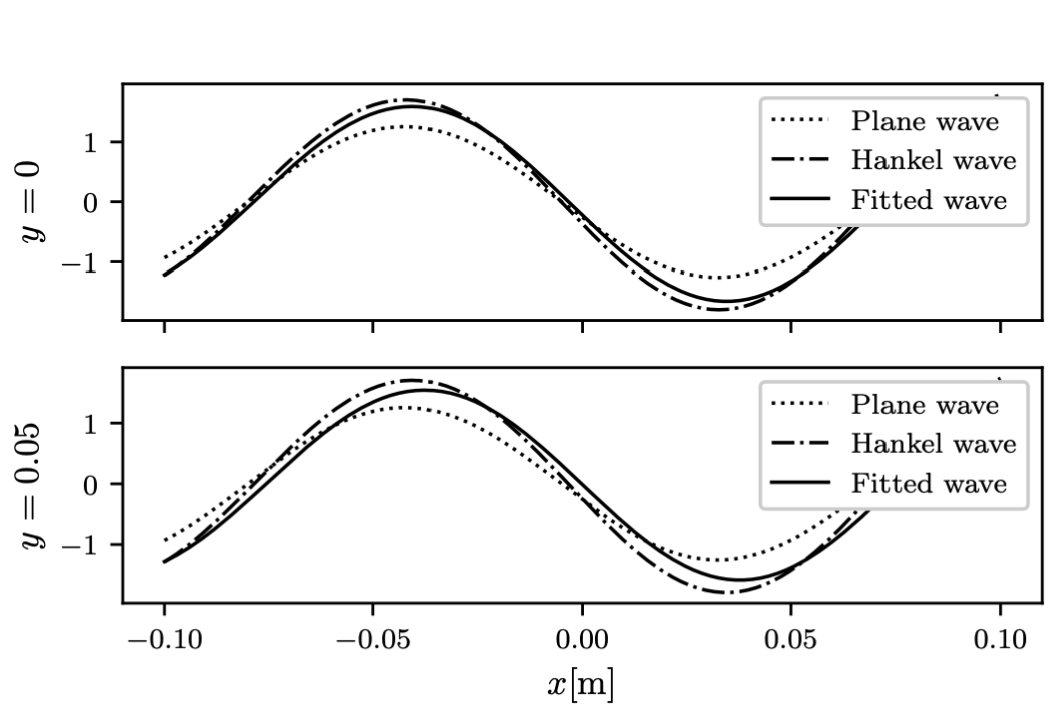}
\caption{Comparison of the three wave models at  $y=0\,\mathrm{m}$ and $y=0.05\,\mathrm{m}$  for the frequency $2\,\mathrm{GHz}.$
}\label{fig:2d_cuts_low}
\end{figure}

\begin{figure}[h!]
\centering
\includegraphics[width=10cm]{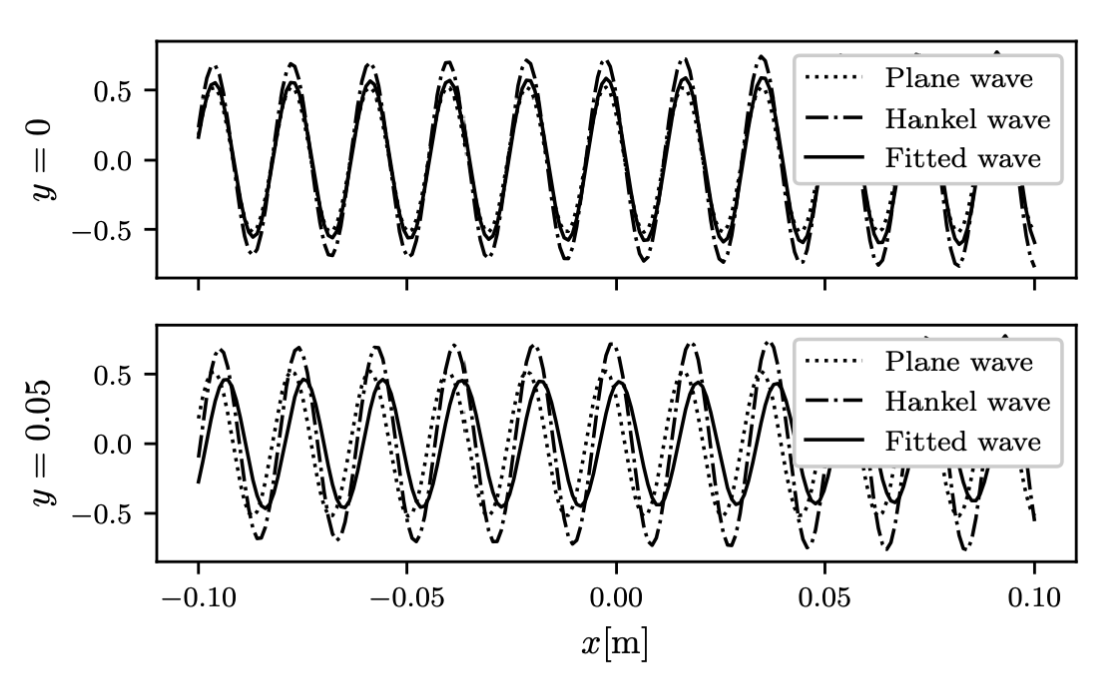}
\caption{ Counterpart of Figure \ref{fig:2d_cuts_low}  for the frequency $16\,\mathrm{GHz}$. }\label{fig:2d_cuts_high}
\end{figure}

\begin{figure}[h!]
\centering
\includegraphics[width=10cm]{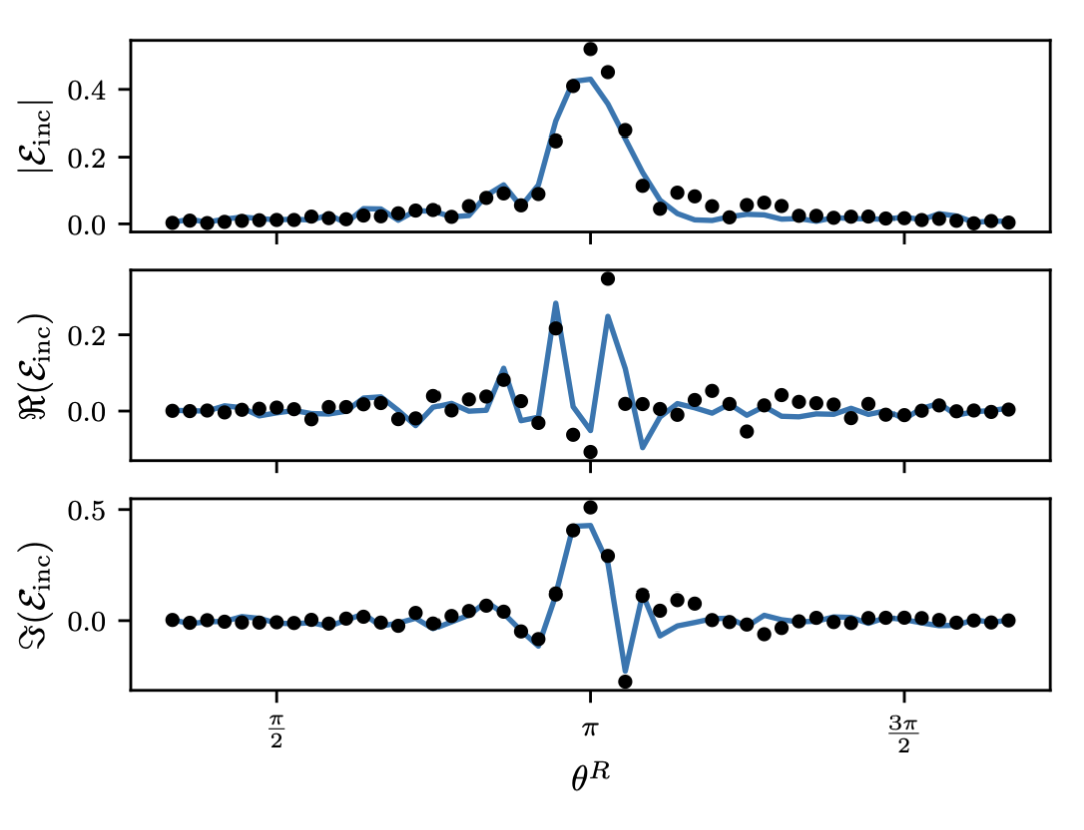}
\caption{Measured incident field for the experiments performed at the highest frequency $16\,\mathrm{GHz}$ (dots) and fitted incident field (line).
Compared to Fig. \ref{fig:3_fitted_low}, the fitting worsens.}
\label{fig:3_fitted_high}
\end{figure}

Figures \ref{fig:2d_comparisons_low}(a)-(c) and \ref{fig:2d_cuts_low} compare the spatial structure of the three proposed fields near the targets: the combination of $14$ modes (\ref{fitted}), the single Hankel wave or isotropic wave defined in (\ref{eq:state_planewave2d}) and the plane wave (\ref{plane}). The plane wave has the smallest amplitude, since we are considering the point furthest from the emitting antenna and the plane wave preserves its amplitude as it advances. The isotropic wave has the largest amplitude, being adjusted to an intermediate point. However the overall shape of the three fields is  similar.

As we increase the frequency of the emitting antenna, the incident wave becomes more anisotropic. Figure \ref{fig:3d_2dwave_low}(b)  represents the fitting for the highest frequency in the data set, $16\,\mathrm{GHz}$.
In this case, the difference between the three approximations is more pronounced, see Figures \ref{fig:2d_comparisons_low}(d)-(f) and \ref{fig:2d_cuts_high}.
Moreover, these fields do not fit so well the experimental data, as shown by Figure \ref{fig:3_fitted_high}.

In the forthcoming sections, incident waves will be fitted by isotropic waves of the form (\ref{eq:state_planewave2d}) since they are slightly simpler to implement. We have already compared with the other two fittings and results are qualitatively the same and omitted for the sake of brevity.

\subsection{Target reconstruction}
\label{sec:td2d}

In the sequel, we will denote by  $\mathcal{E}^\mathrm{inc}_{qkj}$ and $\mathcal{E}^\mathrm{meas}_{qkj}$  the  incident and total electric fields, respectively, measured at $\mathbf{x}^R_j$ when $\mathbf{x}^E_q$ emits at frequency $\nu_k$ in the presence of a target ${\mathcal T}$.
In the same way, we will denote by $\mathcal{E}_{qk}\left(\mathbf{x}^\mathrm{R}_j;\Omega\right)$ the  value of the solution of the forward model for the same experimental conditions but with object $\Omega$  (namely, the solution of (\ref{eq:2ddiel}) or (\ref{eq:2dcond}), depending on the kind of target). As explained in Section \ref{sec:td}, we look for  shapes $\Omega$ such that the difference  between $\mathcal{E}^\mathrm{meas}_{qkj}$ and $\mathcal{E}_{qk}\left(\mathbf{x}^\mathrm{R}_j;\Omega\right)$ is minimized, as
measured by an adequate misfit functional.
The misfit functional for each experiment is defined as  proportional to the $2-$norm distance between the measured data and the synthetic data that would be obtained evaluating at the receivers the solution of the forward problem with objects $\Omega$:
\begin{equation}\label{eq:exp_functional_2D}
\mathcal{J}_{qk}\left({\cal R}\setminus\overline{\Omega}\right)=\frac{1}{2}\sum_{j=1}^{N_\mathrm{meas}}\vert\mathcal{E}_{qk}\left(\mathbf{x}^\mathrm{R}_j;\Omega\right)-\mathcal{E}^\mathrm{meas}_{qkj}
\vert^2.
\end{equation}

When we work with dielectric targets, explicit expressions for the topological derivatives
of such shape functionals with 2D transmission  Helmholtz  problems of the form (\ref{eq:2ddiel}) as constraints are well known  \cite{carpio2008topological,guzina2006misfit}:
\begin{equation}\label{eq:DT_diel2D}
\mathrm{D}_\mathrm{T}\mathcal{J}_{qk}\left(\mathbf{x}\right)=\left(\varepsilon_\mathrm{d}-1\right){\rm Re}\left(\mathcal{U}_{qk}\left(\mathbf{x}\right)\overline{\mathcal{V}_{qk}\left(\mathbf{x}\right)}\right),
\end{equation}
where $\mathcal{U}_{qk}$  solves problem (\ref{eq:2ddiel}) with
$\Omega = \emptyset$, that is,
\begin{equation}\label{eq:state2D}
\left\{
\begin{array}{l}
\Delta\mathcal{U}_{qk}+\kappa_k^2\mathcal{U}_{qk}=0, \qquad \mathrm{in}\; \mathbb R^2,\\
\lim_{\left|\mathbf{x}\right|\to\infty}\sqrt{\left|\mathbf{x}\right|}\left(\frac{\partial\left(\mathcal{U}_{qk}-\mathcal{E}^\mathrm{inc}_{qk}\right)}{\partial \left|\mathbf{x}\right|}-i\kappa_k\left(\mathcal{U}_{qk}-\mathcal{E}^\mathrm{inc}_{qk}\right)\right)=0,
\end{array}
\right.
\end{equation}
and  $\mathcal{V}_{qk}$ is the solution to the following adjoint problem:
\begin{equation}\label{eq:adjoint2D}
\left\{
\begin{array}{ll}
\Delta\mathcal{V}_{qk}+\kappa_k^2\mathcal{V}_{qk}=\displaystyle \sum_{j=1}^{N_{\mathrm{meas}}}(\mathcal{E}_{qkj}^\mathrm{inc}-\mathcal{E}_{qkj}^\mathrm{meas})\delta\left(\mathbf{x}-\mathbf{x}^R_j\right) & \mathrm{in}\; \mathbb R^2,\\
\lim_{\left|\mathbf{x}\right|\to\infty}\sqrt{\left|\mathbf{x}\right|}\left(\frac{\partial\mathcal{V}_{qk}}{\partial \left|\mathbf{x}\right|}+i \kappa_k \mathcal{V}_{qk}\right)=0, &
\end{array}
\right.
\end{equation}
where $\delta$ is Dirac's delta distribution  and $\kappa_k$ is the wavenumber associated to the frequency $\nu_k$, namely $\kappa_k=2\pi\nu_k\sqrt{\mu_0\varepsilon_0}$ (recall that $\mu_0$ and $\varepsilon_0$ are the permeability and the permittivity of vacuum, respectively).
The solution of (\ref{eq:state2D}) is the incident wave itself. Since the Fresnel database only contains values of the incident field at the receivers,  to evaluate the topological derivative (\ref{eq:DT_diel2D}) we need to extend the provided values to the whole plane.
 As already mentioned, we use the approximation (\ref{eq:state_planewave2d}). On the contrary, to solve the adjoint problem  (\ref{eq:adjoint2D}) we only need the values  $\mathcal{E}_{qkj}^\mathrm{meas}$ and $\mathcal{E}_{qkj}^\mathrm{inc}$, which are  stored in the dataset.
The solution of the adjoint problem expressed in terms of fundamental solutions
of the Helmholtz  equation is
\begin{equation}\label{eq:adjoint_form2d}
\mathcal{V}_{qk}\left(\mathbf{x}\right)=\sum\limits_{j=1}^{N_\mathrm{meas}} \left(\mathcal{E}_{qkj}^{\rm inc}-\mathcal{E}^\mathrm{meas}_{qkj}\right)\frac{\imath}{4}H_0^2\left(\kappa_k\left|\mathbf{x}-\mathbf{x}^R_j\right|\right),
\end{equation}
where $H_0^2$ is the Hankel function of second kind and order  zero, that is, an isotropic fundamental solution that radiates from infinity.

For metallic targets, we need the expression of the topological derivative for exterior problems with homogeneous Dirichlet conditions (see \cite{carpio2010determining,samet2003dirichlet}):
\begin{equation}\label{eq:DT_cond2D}
\mathrm{D}_\mathrm{T}\mathcal{J}_{qk}\left(\mathbf{x}\right)={\rm Re}\left(\mathcal{U}_{qk}\left(\mathbf{x}\right)\overline{\mathcal{V}_{qk}\left(\mathbf{x}\right)}\right),
\end{equation}
where $\mathcal{U}_{qk}$ and $\mathcal{V}_{qk}$  solve the same problems (\ref{eq:state2D}) and (\ref{eq:adjoint2D}) as in the case of dielectric targets.   Therefore, $\mathcal{V}_{qk}$  is given by  (\ref{eq:adjoint_form2d}) and $\mathcal{U}_{qk}$ will be  approximated by (\ref{eq:state_planewave2d}).

Notice that the expressions for the topological derivatives for conducting and dielectric targets differ only by a scaling factor  (compare formulas (\ref{eq:DT_diel2D}) and (\ref{eq:DT_cond2D})). Therefore, if we had not known a priori the nature of the objects, we could have used formula  (\ref{eq:DT_cond2D}) to identify the geometry of the targets in both cases.

\begin{figure}[h!]
\centering
\includegraphics[width=10cm]{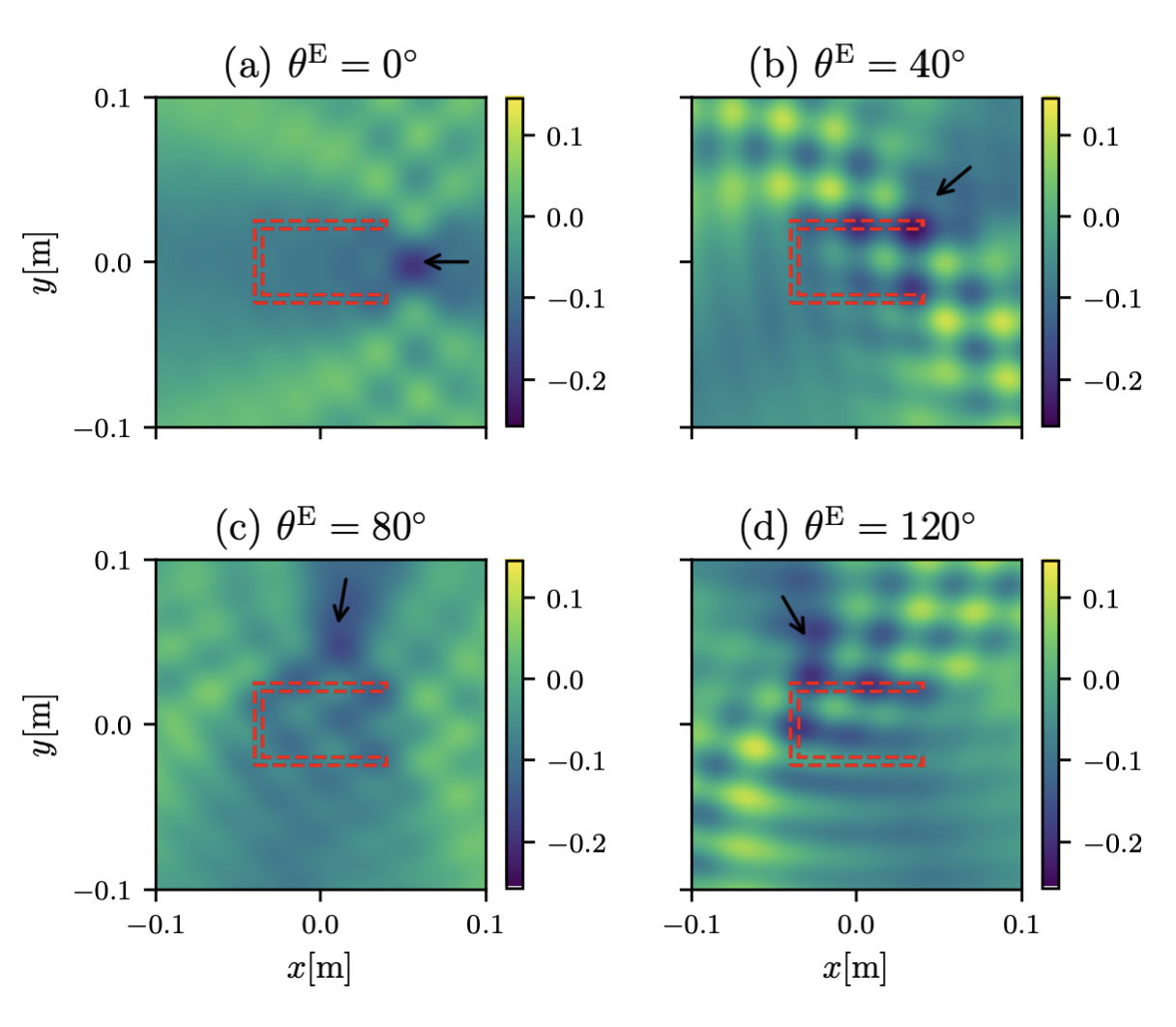}
\caption{Topological derivative (\ref{eq:DT_cond2D}) for $4$ different antenna positions described by the angle $\theta^E$ in (\ref{position}) at frequency $6\,\mathrm{GHz}$. Red contours show the true target.}
\label{fig:U_single_exps2}
\end{figure}

Figure \ref{fig:U_single_exps2} shows the topological derivative for $4$ experiments with  the emitting antenna located at different positions  described by the indicated angles $\theta^E$ (see (\ref{position})). We represent the vector joining the  antenna position to the origin of coordinates by an arrow. The minima of the derivative coincide sometimes with a point of the object. We choose the U-shaped target  because it is the most  demanding shape, extremely non-convex.

To use all the available information, we consider misfit functionals combining different experiments and different frequencies. We denote by
\begin{equation} \label{multidirection_td}
\mathcal{J}_{k}\left({\cal R}\setminus\overline{\Omega}\right)=\frac{1}{N_\mathrm{exp}}\sum_{q=1}^{N_\mathrm{exp}}\mathcal{J}_{qk}\left({\cal R}\setminus\overline{\Omega}\right)
\end{equation}
the single-frequency functional combining all the rotations of the target, that is, the $N_\mathrm{exp}$ experiments corresponding to the $k$-th frequency.  By linearity, it is clear that the topological derivative of $\mathcal{J}_{k}$ is nothing but the linear combination of the individual ones:
\begin{equation}\label{eq:DTmulti}
\mathrm{D}_\mathrm{T}\mathcal{J}_k=\frac{1}{N_{\rm exp}}\,\sum_{q=1}^{N_{\rm exp}} \mathrm{D}_\mathrm{T}\mathcal{J}_{qk}.
\end{equation}

\begin{figure}[h!]
\centering
\includegraphics[width=10cm]{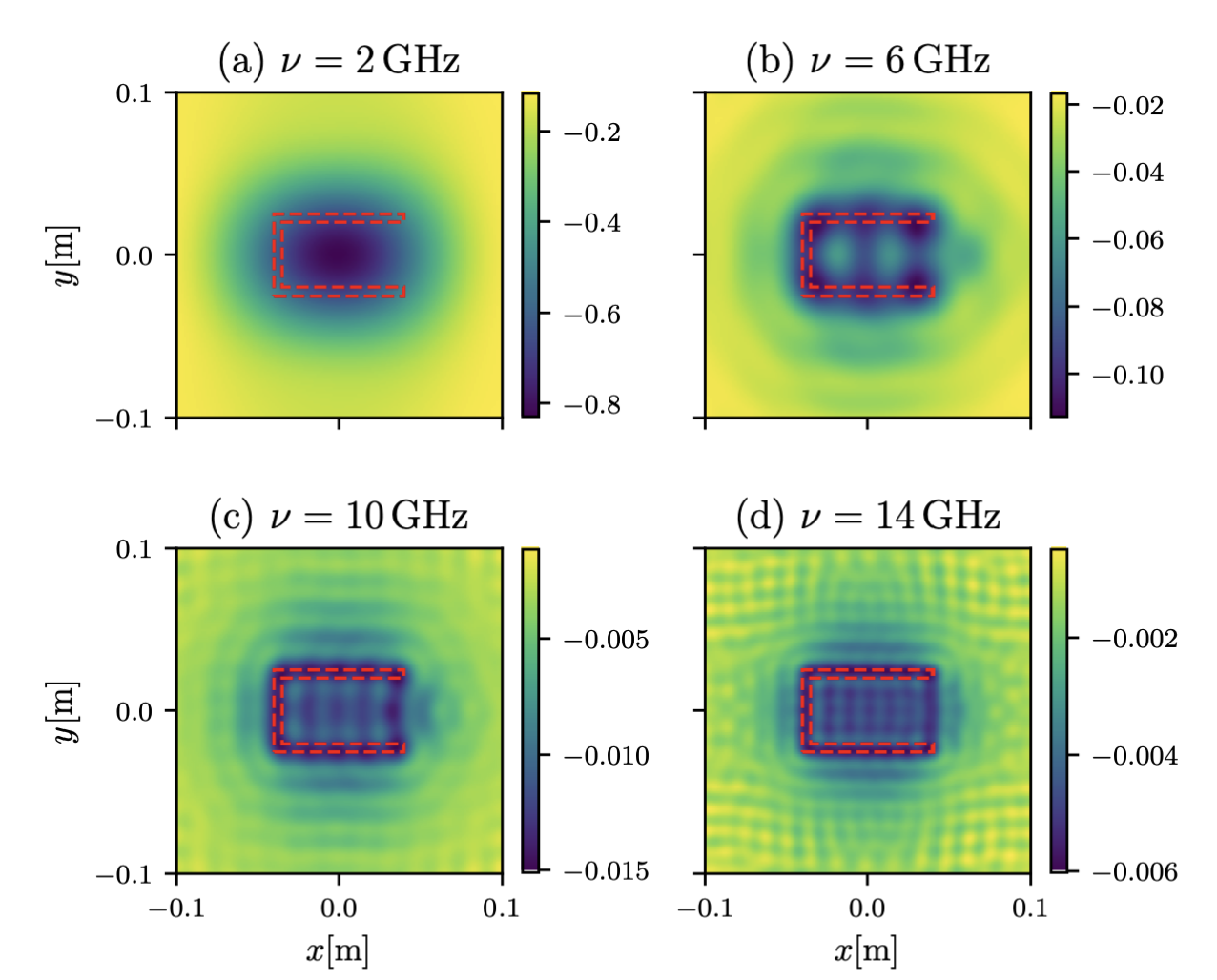}
\caption{Single frequency topological derivative  (\ref{eq:DTmulti})  for $4$ different frequencies. Red contours mark the true object. }
\label{fig:U_single_freq}
\end{figure}

Figure \ref{fig:U_single_freq} shows that combining information from all the positions of the emitting antenna we improve the results since the largest negative values (dark blue colors in the plots) somehow identify the target.  For low frequencies we approximate of the size and position of the object but obtain little information about the shape. For high frequencies we recover some features of the shape, but spurious minima do appear.

To average information from different frequencies $\nu_k,\ k=1,\dots, N_{\rm freq}$, we consider the following linear combination of the single-frequency topological derivatives $\mathrm{D}_\mathrm{T}\mathcal{J}_k$ defined in (\ref{eq:DTmulti}):
\begin{equation} \label{multifrequency_td}
\mathrm{D}_{\mathrm{T}}\mathcal{J}\left(\mathbf{x}\right) = \frac{1}{N_\mathrm{freq}}\sum_{k=1}^{N_\mathrm{freq}}\frac{\mathrm{D}_{\mathrm{T}}\mathcal{J}_{k}\left(\mathbf{x}\right)}{\vert\min_{\mathbf{y}\in {\cal R}^{\rm insp}} \mathrm{D}_{\mathrm{T}}\mathcal{J}_{k}\left(\mathbf{y}\right)\vert},
\end{equation}
where the region ${\cal R}^{\rm insp}$ that appears in the weights $\vert\min_{\mathbf{y}\in{\cal R}^{\rm insp}}\mathrm{D}_{\mathrm{T}}\mathcal{J}_{k}\left(\mathbf{y}\right)\vert^{-1}$ is the inspection region, namely, the region  where the topological derivative will be evaluated (see Figure \ref{fig:2d_diagram}). The choice of this kind of weights was firstly proposed in \cite{funes2016defect} and it is motivated by the fact that although when considering different directions we expect the total energy/amplitude of the incident and total waves  to be similar  (as can be observed in Figure \ref{fig:U_single_exps2}), this is not true in general for different frequencies,  as can be observed in Figure \ref{fig:U_single_freq}, where it ranges from $-0.2$ to $0.1$ for $\nu=2$ GHz, while the range reduces to $[-0.006,0.003]$ for $\nu=14$ GHz.

\begin{figure}[h!]
\centering
\includegraphics[width=10cm]{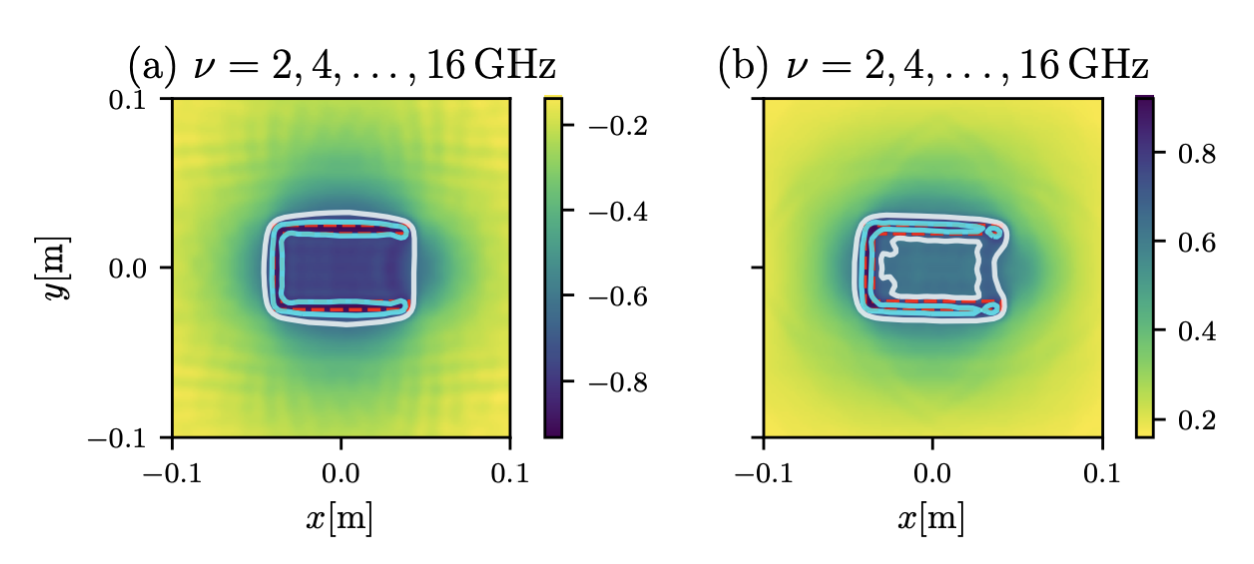}
\caption{Multi-frequency (a) topological derivative (\ref{multifrequency_td}) and (b)  topological energy (\ref{multifrequency_et})  for the U-shaped target.
The true targets are represented by dashed red contours.   Solid white and cyan
curves represent contour levels for $\lambda=0.7$ and $\lambda=0.9$, respectively.
}
\label{fig:U_multi_freq}
\end{figure}

Figure \ref{fig:U_multi_freq}(a) superimposes information from all the frequencies considered in the dataset  {\tt uTM\_shaped.exp} (name of the file in \cite{belkebir2001testing},  see also Table \ref{tab:exps}).  We mark by a white line the boundary of the approximate domain $\Omega_\lambda$ defined in (\ref{eq:dom_app}) for $\lambda = 0.7$, while the cyan curve corresponds to $\lambda=0.9$.   This example shows that  for this particular target,  the topological derivative is very good at capturing the convex hull of the shape for a wide range of values of $\lambda$, and that for some of such values we do obtain rather accurate approximations of the correct shape of the target.

\begin{figure}[h!]
\centering
\includegraphics[width=10cm]{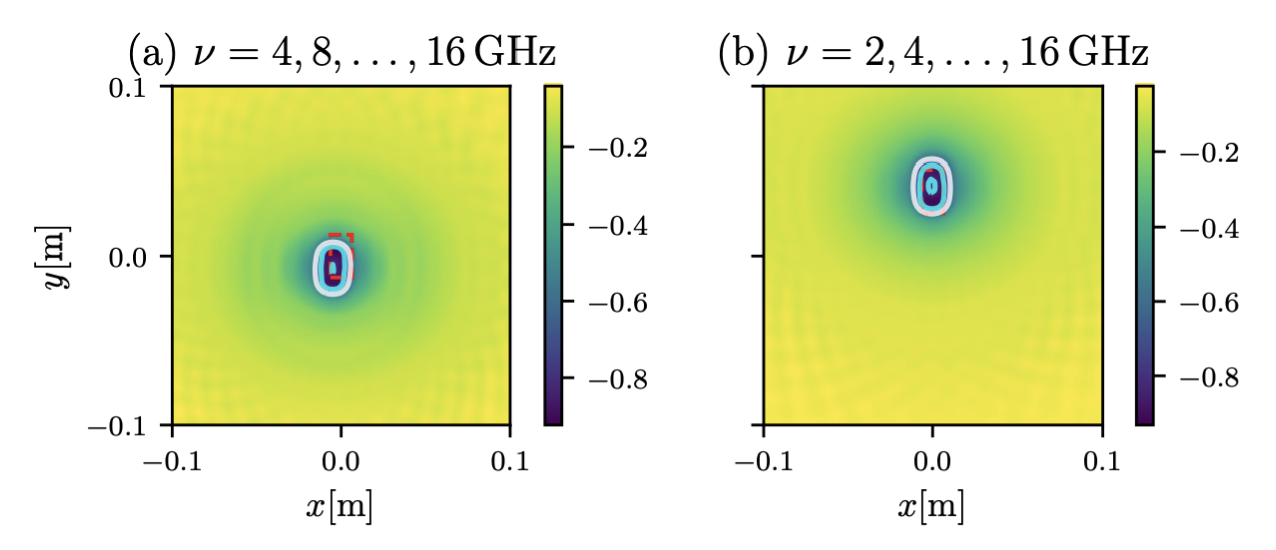}
\caption{Multi-frequency topological derivative (\ref{multifrequency_td})
for (a) the centered rectangular target and (b) the off-centered rectangular target.
The true targets are represented by dashed red contours,
while solid white  and cyan contours represent the levels $\lambda =0.7$ and $\lambda=0.9$, repectively.
}
 \label{fig:RectCent_multi_freq}
\end{figure}

Figure \ref{fig:RectCent_multi_freq} depicts the multi-frequency topological derivatives corresponding to the other two conducting targets (datasets  {\tt rectTM\_cent.exp} and {\tt rectTM\_dece.exp} respectively). Given the accuracy in  the reconstruction of both the U-shaped target as well as the off-centered rectangle, it seems that the mismatch observed for the centered rectangle must be due to an experimental error. Indeed, other papers processing the same data with different techniques  made the same observation (compare, for instance, with \cite{bloemenkamp2001inversion,duchene2001inversion,marklein2001linear}).

\begin{figure}[h!]
\centering
\includegraphics[width=10cm]{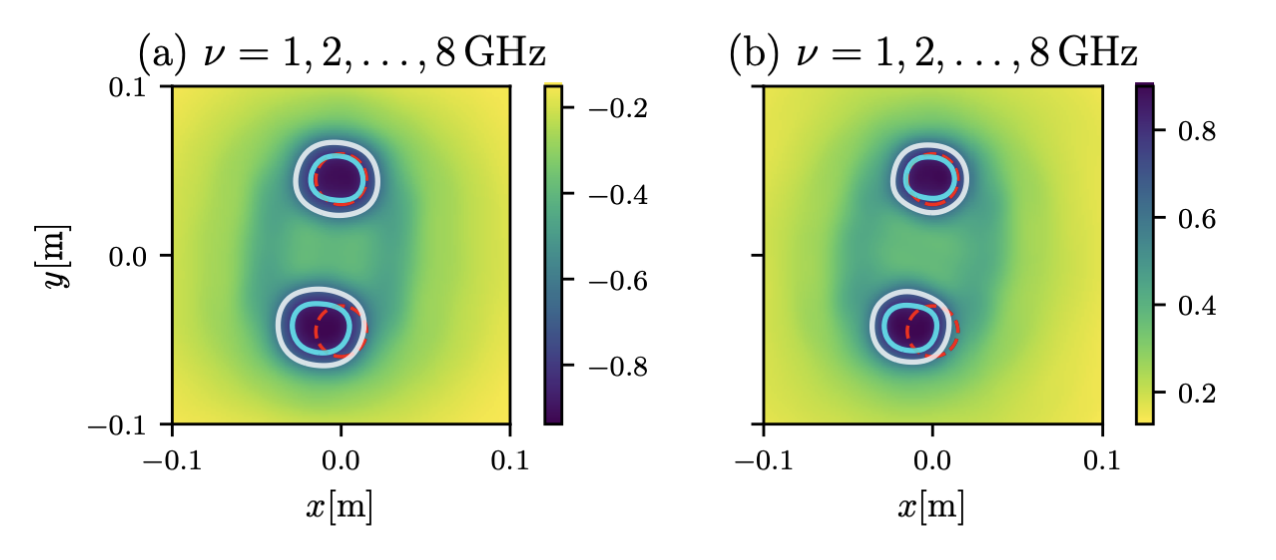}
\caption{Multi-frequency (a) topological derivative (\ref{multifrequency_td})
and (b) topological energy (\ref{multifrequency_et}) for the two
dielectric cylinders. The true targets are represented by dashed
red contours, while the solid white  and cyan contours represent the levels $\lambda =0.7$ and $\lambda=0.9$, respectively.}
\label{fig:2c8f_multi_freq}
\end{figure}

We encounter a similar phenomenon for dataset {\tt twodielTM\_8f.exp}, corresponding to the target formed by two dielectric cylinders. Figure \ref{fig:2c8f_multi_freq}(a) shows that the multi-frequency topological derivative correctly finds the shape, size and number of objects. However there is a mismatch in their position/orientation. This discrepancy reminds of a solid translation and rotation, so it may related to experimental errors/noise or aberrations of the imaging setup. This fact was previously observed by other research groups (see \cite{bloemenkamp2001inversion} or \cite{marklein2001linear} for instance).
An additional dataset corresponding to the same target with a different set of frequencies {\tt twodielTM\_4f.exp} produced no meaningful reconstructions. Similar problems were encountered in \cite{bloemenkamp2001inversion,testorf2001imaging} when analyzing the same dataset. As mentioned in \cite{testorf2001imaging}, some magnitudes are more than double than the rest of the dataset, which could indicate some kind of data corruption.

\begin{figure}[h!]
\centering
\includegraphics[width=10cm]{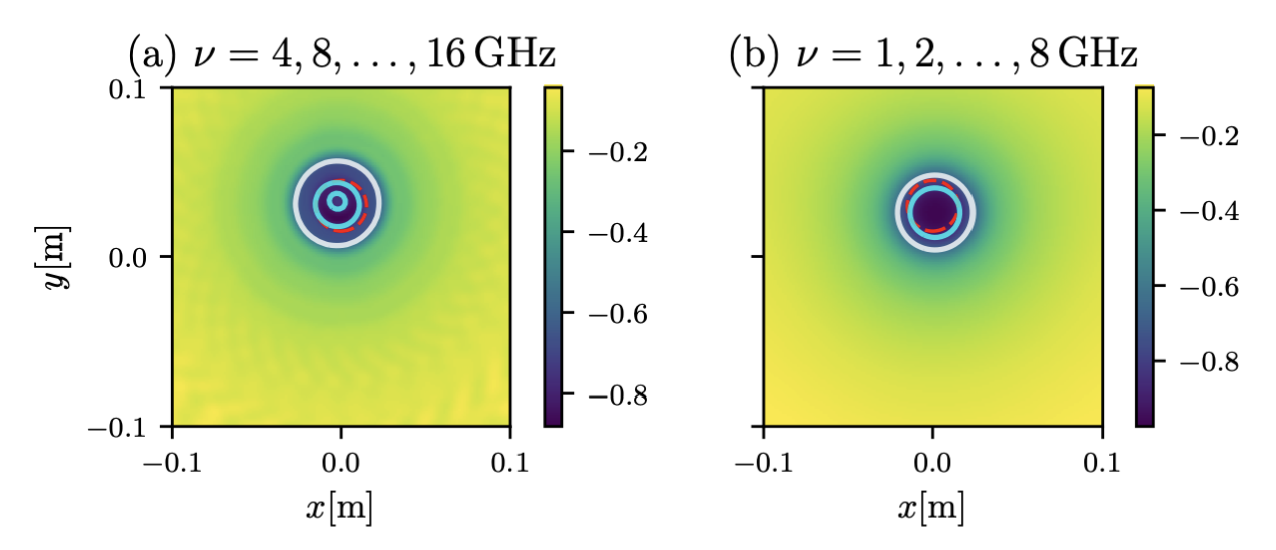}
\caption{Multi-frequency topological derivatives (\ref{multifrequency_td})
for the dielectric cylinder and different frequency choices.
The true targets are represented by dashed
red contours, while the solid white and cyan contours represent the levels $\lambda =0.7$ and $\lambda=0.9$, repectively.}
 \label{fig:1c4f_multi_freq}
\end{figure}

Finally, we have analyzed datasets {\tt dielTM\_dec4f.exp} and {\tt dielTM\_dec8f.exp}, see Figure \ref{fig:1c4f_multi_freq}. The reconstruction improves
varying the frequencies, since (a) presents a hole for $\lambda=0.9$ (cyan countour) and overestimates the size of the cylinder for $\lambda=0.7$  (white contour).

All the results presented were computed on a $4$-cores desktop computer using an $N\times N$ grid with $N=100$. We processed the datasets with $4$ frequencies in $6$ seconds whereas the datasets with $8$ frequencies run in $12$ seconds. Note that this is a one-step method which uses closed-form expressions for topological derivatives and auxiliary fields. As a result, the computation time scales with the number of nodes $M=N^2$, that is, the time complexity is $O\left(M\right)$. The value of $N=100$ was chosen only for visualization reasons, to be able to plot the fields with enough resolution, but lower values would be acceptable.

A closely related but somehow simpler indicator function,  the topological energy,  only uses the norms of the adjoint (\ref{eq:adjoint2D}) and forward fields (\ref{eq:state2D}) (see \cite{dominguez2010non}):
\begin{equation}
\mathrm{E}_\mathrm{T}\mathcal{J}_{qk}(\mathbf{x})= \left|\mathcal{U}_{qk}\left(\mathbf{x}\right)\right|^2\left|\mathcal{V}_{qk}\left(\mathbf{x}\right)\right|^2.
\end{equation}
It is useful when the topological derivative displays an oscillatory behavior which is difficult to interpret (which is not really the case here - it is for high frequencies, see figure \ref{fig:U_single_freq}(c) and \ref{fig:U_single_freq}(d)), and also when the nature of the targets  is not know beforehand or when targets of different nature are simultaneously present in the media (see \cite{rapun2020}). When the incident field is approximated by a plane wave, we have $ \left|\mathcal{U}_{qk}\left(\mathbf{x}\right)\right|^2=\mathrm{constant} $, which decreases the computational cost (and allows us to skip the whole fitting of $\mathcal{U}_\mathrm{inc}$ process). The multi-frequency version would be
\begin{eqnarray}
\mathrm{E}_\mathrm{T}\mathcal{J}(\mathbf{x}) = \frac{1}{N_\mathrm{freq}}\sum_{k=1}^{N_\mathrm{freq}}\frac{
\mathrm{E}_\mathrm{T}\mathcal{J}_{k}(\mathbf{x})}
{{\max}_{ \mathbf{y}\in {\cal R}^{\rm insp} }\mathrm{E}_{\mathrm{T}}\mathcal{J}_{k}(\mathbf{y})},
\quad
\mathrm{E}_\mathrm{T}\mathcal{J}_{k}=\frac{1}{N_\mathrm{exp}}\sum_{q=1}^{N_\mathrm{exp}}\mathrm{E}_\mathrm{T}\mathcal{J}_{qk}.
\label{multifrequency_et}
\end{eqnarray}
In Figures   \ref{fig:U_multi_freq} and \ref{fig:2c8f_multi_freq} we compare results obtained with the multi-frequency topological energies and derivatives.  For the case of the topological energy, the approximated domain is defined by the counterpart  to (\ref{eq:dom_app}), namely,
\begin{equation}\label{eq:omegaap2}
\Omega_{\mathrm{app},\lambda}:=\left\lbrace \mathbf{x}\in {\cal R}^{\rm insp} \;\mathrm{s.t.}\;\mathrm{E}_\mathrm{T}\mathcal{J}\left(\mathbf{x}\right)\ge \lambda\max_{\mathbf{y}\in {\cal R}^{\rm insp}}\mathrm{D}_\mathrm{T}\mathcal{J}\left(\mathbf{y}\right) \right\rbrace.
\end{equation}
For an easier comparison between topological derivatives and energies, we have reversed the colormap scale for the second ones, namely, dark blue color corresponds to the minimum for the topological derivative plots while it corresponds to the maximum for the topological energy ones.

The maxima of the topological energies are somewhat more pronounced than the minima of the topological derivatives, which makes sense since the topological energy scales as the square of the topological derivative. Indeed, the reconstructions for the two values of $\lambda$ keep a similar shape for the U-shaped target, and remain closer for the two cylinders, meaning that they are more robust with respect to $\lambda$.

\section{Application to the three dimensional Fresnel database}
\label{sec:3D}

\begin{figure}[h!]
\centering
\includegraphics[scale=0.9]{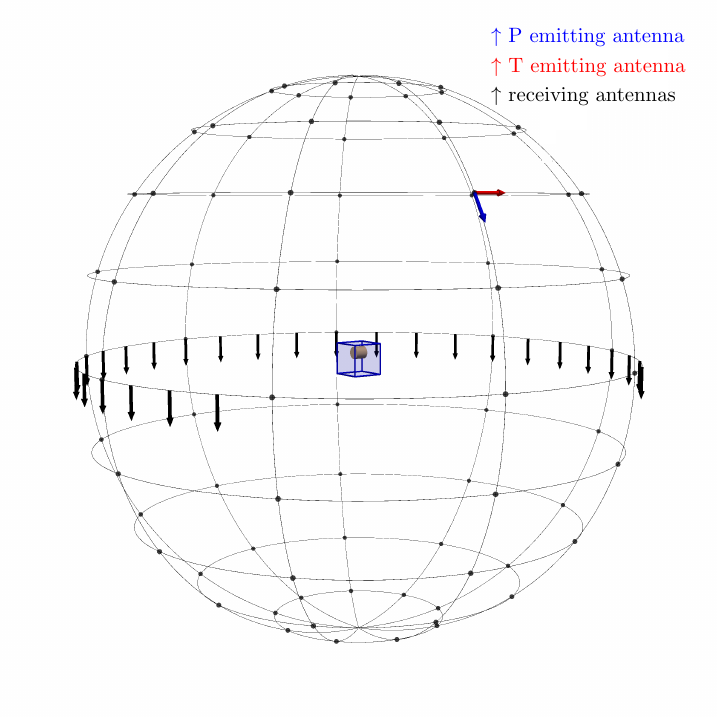}
\caption{Scaled representation of the location of the receiving and emitting antennas for a generic experiment, as well as the size of the inspection zone compared with the size of the objects. The black dots represent the possible positions for the emitting antenna. The receiving antennas will move according to the emitting antenna position.
}
\label{fig:3d_diagram}
\end{figure}

The 3D Fresnel database \cite{geffrin2009continuing} comprises a much more general geometric configuration for the experiments.
In this database, targets are irradiated from a given point on a sphere of radius $R^\mathrm{E} = 1.796 \,\mathrm{m}$. The incident and total fields are measured at $N_\mathrm{meas} = 27$ receivers located on the horizontal plane, see Figure \ref{fig:3d_diagram}.

To represent the position and emission directions of the antennas, we use a spherical coordinate system $\left(\mathrm{r},\theta,\phi\right)$, where $0\le\mathrm{r}<\infty$ is the distance to the origin, $0^\circ\le\theta<360^\circ$ is the azimuth angle, that is, the angle in the horizontal plane, and $0^\circ\le\phi\le 180^\circ$ is the altitude angle measured with respect to the vertical positive direction:
\begin{eqnarray*}
\mathrm{r} =\sqrt{x^2 + y^2 + z^2}, \quad
\theta = \arctan\left(\frac{y}{x}\right), \quad
\phi = \arctan\left(\frac{\sqrt{x^2 + y^2}}{z}\right).
\end{eqnarray*}
Associated to this coordinate system is the set of unit vectors:
\begin{eqnarray*}
\mathbf{u}_\mathrm{r}(\theta,\phi) =\cos\theta\sin\phi\,\mathbf{i}+\sin\theta\sin\phi\,\mathbf{j}+\cos\phi\,\mathbf{k},
\\
\mathbf{u}_\theta(\theta,\phi) = -\sin\theta\,\mathbf{i}+\cos\theta\,\mathbf{j},
\\
\mathbf{u}_\phi (\theta,\phi) = \cos\theta\cos\phi\,\mathbf{i}+\sin\theta\cos\phi\,\mathbf{j}-\sin\phi\,\mathbf{k}.
\end{eqnarray*}
As explained in \cite{geffrin2009continuing}, the emitting antenna can be situated at $N_{\phi^\mathrm{E}}=9$ different altitude angles $\phi^\mathrm{E}$, ranging from $\phi^\mathrm{E}=18^{\circ}$ to $\phi^\mathrm{E}=162^{\circ}$ at increments of $18^{\circ}$, and at $N_{\theta^\mathrm{E}} = 9$ different azimuth angles $\theta^\mathrm{E}$, ranging from $\theta^\mathrm{E}=40^{\circ}$ to $\theta^\mathrm{E}=360^{\circ}$ at increments of $40^{\circ}$. We denote these positions by $\mathbf{x}^\mathrm{E}_{pq}$
defined as
\[
\mathbf{x}^\mathrm{E}_{pq} = R^\mathrm{E}\mathbf{u}_\mathrm{r}\left(\theta^\mathrm{E}_p,\phi^\mathrm{E}_q\right)
\]
with $\theta^\mathrm{E}_p$ and $\phi^\mathrm{E}_q$   indexed as:
\[
\theta^\mathrm{E}_p = 40^{\circ}p,\quad p =1,\dots,N_{\theta^\mathrm{E}},
\qquad
\phi^\mathrm{E}_q = 18^{\circ}q,\quad q =1,\dots,N_{\phi^\mathrm{E}}.
\]
These $9\times 9=81$ positions are represented in Figure \ref{fig:3d_diagram} by small black spheres.

Two types of polarizations are considered for the emitting antennas. The so-called parallel polarization (PP), where the polarization vector of the emitting antenna is equal to the unit vector along the altitude direction:
\[
\mathbf{p}^\mathrm{E}=\mathbf{u}_\phi ,
\]
and the transverse polarization (TP), where the polarization vector is equal to the unit vector along the azimuth direction:
\[
\mathbf{p}^\mathrm{E}=\mathbf{u}_\theta.
\]
An example of both polarizations is represented in Figure \ref{fig:3d_diagram}, where a red arrow is used for the transverse polarization and a blue one for the parallel polarization.

For each position of the emitting antenna, the measuring antennas are situated in the horizontal plane (hence $\phi_\mathrm{R}=90^{\circ}$), on a circumference of radius $R^\mathrm{R}=R^\mathrm{E}$ and at $N_{\theta^\mathrm{R}}=27$ azimuth angles $\theta^\mathrm{R}$ varying from \mbox{$\theta^\mathrm{R}=\theta^\mathrm{E}+50^{\circ}$} to $\theta_{R}=\theta^\mathrm{E}+310^{\circ}$ at increments of $10^{\circ}$. We will denote these positions by  $\mathbf{x}^\mathrm{R}_{j}$, such that:
\[
\mathbf{x}^\mathrm{R}_{j} = R^\mathrm{R}\mathbf{u}_\mathrm{r}\left(\theta^\mathrm{R}_j,90^\circ\right),
\]
with $\theta^\mathrm{R}_j$ being indexed as:
\[
\theta^\mathrm{R}_j = \theta^\mathrm{E}+50^{\circ}j,\quad j =1,\dots,N_{\theta^\mathrm{R}},
\]
\noindent that is, the position of the receiving antennas with respect to the azimuth angle of the emitting antenna is fixed. All the receiving antennas are polarized in the parallel direction, and since they are in the horizontal plane we have:
\[
\mathbf{p}^\mathrm{R}=-\mathbf{k}.
\]
The receiving antennas are represented in Figure \ref{fig:3d_diagram} by black arrows.

Finally, for each emitting position and polarization, $N_\mathrm{freq}=21$ frequencies were used, ranging from $3\,\mathrm{GHz}$ to $8\,\mathrm{GHz}$ with a $0.25\,\mathrm{GHz}$ increment step. We will index these frequencies as:
\[
\nu_k = 3+0.25\left(k-1\right)\,\mathrm{GHz},\quad k=1,\dots,N_\mathrm{freq},
\]
and the corresponding wave-numbers in  vacuum as:
\[
\kappa_k = 2\pi \nu_k \sqrt{\mu_0\varepsilon_0},\quad k=1,\dots,N_\mathrm{freq}.
\]

\begin{figure}[h!]
\centering
\includegraphics[width=10cm]{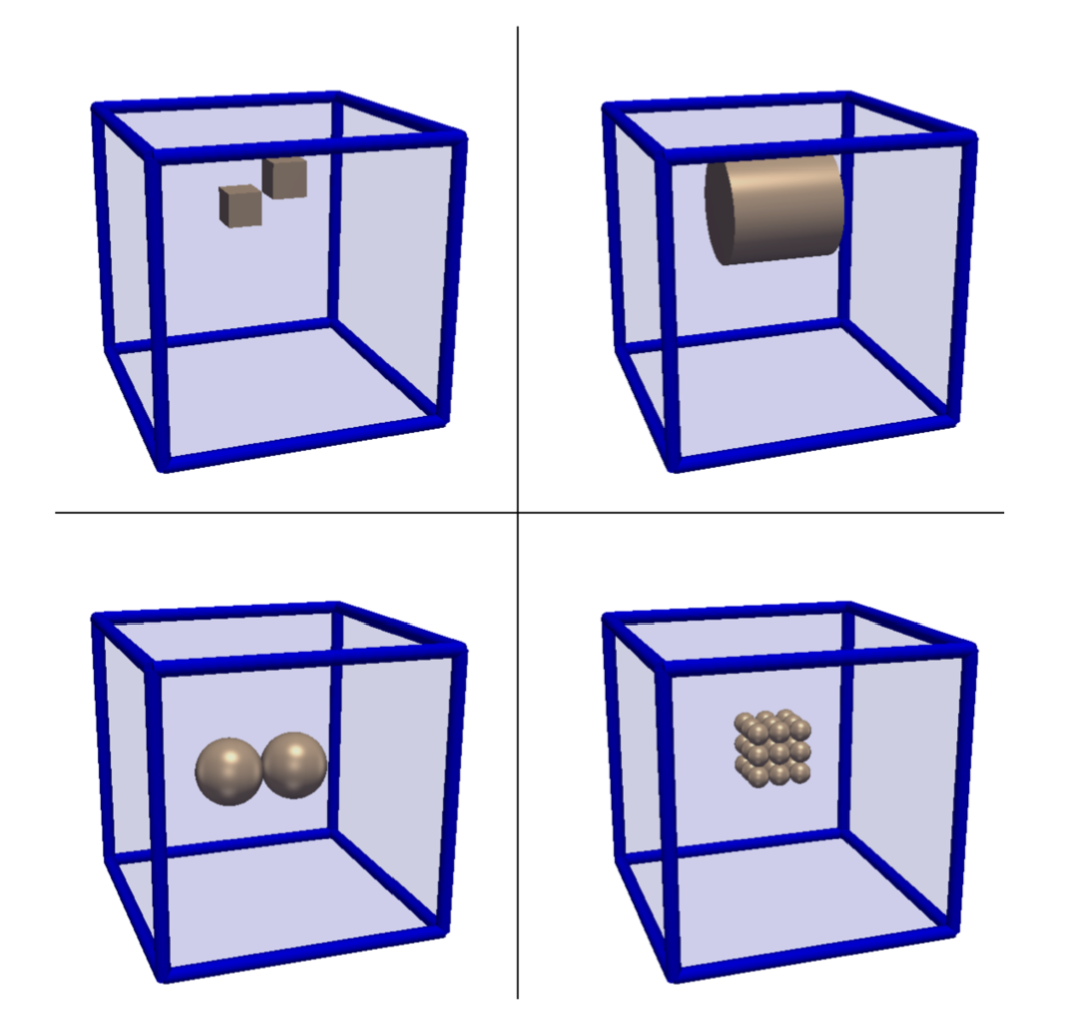}
\caption{ 3D targets represented along with the fixed inspection zone. 
}
\label{fig:3d_targets}
\end{figure}

Five different dielectric (with relative magnetic permeability $\mu_\mathrm{r}=1$) targets were studied. Two spheres with relative electric permittivity of $\varepsilon_\mathrm{r}=2.6$, two cubes with $\varepsilon_\mathrm{r}=2.3$, a cylinder with $\varepsilon_\mathrm{r}=3.05$, a cubic array of nine spheres with $\varepsilon_\mathrm{r}=2.6$, and finally a mysterious target whose shape was not known beforehand. The first $4$ targets are shown in figure \ref{fig:3d_targets}, where the inspection zone ${\cal R}^{\rm insp}=[-0.1,0.1]^3$
 is also represented. The same inspection zone is represented in Figure \ref{fig:3d_diagram} to compare the actual size of the targets with the distance of the emitting and receiving antennas.

Taking into account the general geometrical configuration of the experiment as well as the targets,  the three components of the electrical field in the Maxwell equations are coupled,  and the full 3D Maxwell models (\ref{eq:3ddiel}) and (\ref{eq:3dcond}) need to be used.
To implement the imaging strategy described in Section \ref{sec:td} we need to
introduce the proper topological fields.

\subsection{Topological fields}
\label{sec:td3d}

In the sequel, we will denote by  $\boldsymbol{\mathcal{E}}^\mathrm{inc}_{pqkj}$ and $\boldsymbol{\mathcal{E}}^\mathrm{meas}_{pqkj}$  the measurement of  the incident and the total electric fields  performed at position  $\mathbf{x}^\mathrm{R}_{j}$  when irradiating the target from position $\mathbf{x}^\mathrm{E}_{pq}$ at frequency  $\nu_k$. In the same way, we will denote by  $\boldsymbol{\mathcal{E}}_{pqk}\left(\mathbf{x}^\mathrm{R}_j;\Omega\right)$ the solution of (\ref{eq:3ddiel}) evaluated at the receiver
$\mathbf{x}^\mathrm{R}_j$ for the same experimental conditions but with a generic object $\Omega$.

The misfit functional for each experiment is defined as one half of the  $2-$norm  distance between the measured data and the synthetic data that would be obtained evaluating at the receivers the solution of the forward problem with object $\Omega$. Note that we do not know the full three dimensional electric  vector field at the measuring points, but only its component along $\mathbf{k}$. The misfit functional   $\mathcal{J}_{pqk}$  is then
\begin{equation}\label{eq:3D_misfit}
\mathcal{J}_{pqk}\left({\cal R}\setminus\overline{\Omega}\right)=\frac{1}{2}\sum_{j=1}^{N_{\theta^\mathrm{R}}}\vert\left(\boldsymbol{\mathcal{E}}_{pqk}\left(\mathbf{x}^\mathrm{R}_j;\Omega\right)-\boldsymbol{\mathcal{E}}^\mathrm{meas}_{pqkj}\right)\cdot\mathbf{k}
\vert^2.
\end{equation}

Following the derivations in \cite{le2017topological,masmoudi2005maxwell}, we find that the topological derivative of the functional (\ref{eq:3D_misfit}) is:
\begin{equation}\label{eq:DTmaxwell}
\mathrm{D}_\mathrm{T}\mathcal{J}_{pqk}\left(\mathbf{x}\right) = -3{\rm Re}\left(\kappa_{k}^{2}\frac{\left(\varepsilon_\mathrm{r}-1\right)}{\left(\varepsilon_{\mathrm{r}}+2\right)}\,\boldsymbol{\mathcal{U}}_{pqk}\left(\mathbf{x}\right)\cdot\overline{\boldsymbol{\mathcal{V}}_{pqk}\left(\mathbf{x}\right)}\right),
\end{equation}
where  $\boldsymbol{\mathcal{U}}_{pqk}$ is the solution of
\begin{equation}\label{eq:maxwelstate}
\left\{
\begin{array}{l}
\mathbf{curl}\,\mathbf{curl}\,\boldsymbol{\mathcal{U}}_{pqk}-\kappa_k^2\boldsymbol{\mathcal{U}}_{pqk}=\mathbf{0},\qquad\mathrm{in}\; \mathbb R^3,\\
\lim_{\left|x\right|\to\infty}\left|\mathbf{x}\right|\left|\mathbf{curl}\left(\boldsymbol{\mathcal{U}}_{pqk}-\boldsymbol{\mathcal{E}}^\mathrm{inc}_{pqk}\right)\times\frac{\mathbf{x}}{\left|\mathbf{x}\right|}-i\kappa_k\left(\boldsymbol{\mathcal{U}}_{pqk}-\boldsymbol{\mathcal{E}}^\mathrm{inc}_{pqk}\right)\right|=0,
\end{array}
\right.
\end{equation}
and  $\boldsymbol{\mathcal{V}}_{pqk}$  in (\ref{eq:DTmaxwell}) is the adjoint field, and it solves:
\begin{equation}\label{eq:maxweladjoint}
\left\{
\begin{array}{l}
\mathbf{curl}\,\mathbf{curl}\,\boldsymbol{\mathcal{V}}_{pqk}-\kappa_k^2\boldsymbol{\mathcal{V}}_{pqk}=\sum\limits_{j=1}^{N_{\theta^\mathrm{R}}}(\boldsymbol{\mathcal{E}}_{pqkj}^\mathrm{inc}-\boldsymbol{\mathcal{E}}_{pqkj}^\mathrm{meas})\cdot\mathbf{k}\,\delta\left(\mathbf{x}-\mathbf{x}^\mathrm{R}_j\right), \ \mathrm{in}\;\mathbb R^3\\
\lim_{\left|x\right|\to\infty}\left|\mathbf{x}\right|\left|\mathbf{curl}\,\boldsymbol{\mathcal{V}}_{pqk}\times\frac{\mathbf{x}}{\left|\mathbf{x}\right|}+i\kappa_k\boldsymbol{\mathcal{V}}_{pqk}\right|=0,
\end{array}
\right.
\end{equation}

Problem (\ref{eq:maxwelstate}) is exactly (\ref{eq:3ddiel}) without a target, that is, $\Omega=\emptyset.$
The solution is the incident wave, which we only know at the measuring locations.
However, the authors of the database performed a calibration in the measurements in such a way that the incident field can be approximated by a plane wave with unitary amplitude and zero phase at the origin. This way, for the parallel polarization case
\[
\boldsymbol{\mathcal{U}}_{pqk}\left(\mathbf{x}\right) =\mathbf{u}_{\phi}\left(\theta_p,\phi_q\right)e^{-i \kappa_k\mathbf{x}\cdot\mathbf{u}_\mathrm{r}\left(\theta_p,\phi_q\right)},
\]
whereas for the transverse polarization case
\[
\boldsymbol{\mathcal{U}}_{pqk}\left(\mathbf{x}\right) =\mathbf{u}_{\theta}\left(\theta_p,\phi_q\right)e^{-i \kappa_k\mathbf{x}\cdot\mathbf{u}_\mathrm{r}\left(\theta_p,\phi_q\right)}.
\]

As it happened with the 2D database, the fitting or approximation performed to obtain an expression for  $\boldsymbol{\mathcal{U}}_{pqk}$  does not play a role in the computation of the adjoint field, which only  depends on  the measurements of the incident field at the receptors provided in the database. The solution to the adjoint problem (\ref{eq:maxweladjoint}) is given by
\[
\boldsymbol{\mathcal{V}}_{pqk}\left(\mathbf{x}\right)=\frac{\imath}{4\kappa_k^2}\sum_{j=1}^{N_{\theta^\mathrm{R}}}\mathbf{curl}\,\mathbf{curl}\left(H_0^2\left(\kappa_k\vert\mathbf{x}-\mathbf{x}^\mathrm{R}_j\vert\right)\left(\boldsymbol{\mathcal{E}}_{pqkj}^\mathrm{inc}-\boldsymbol{\mathcal{E}}_{pqkj}^\mathrm{meas}\right)\cdot\mathbf{k}\right).
\]

As we did for the 2D case, we compute the single-frequency topological derivative as the average of the topological derivative for each emission position, that is
\begin{eqnarray}
\mathrm{D}_\mathrm{T}\mathcal{J}_{k}\left(\mathbf{x}\right) = \frac{1}{N_{\theta^\mathrm{E}}N_{\phi^\mathrm{E}}}\sum	_{p=1}^{N_{\theta^\mathrm{E}}}\sum_{q=1}^{N_{\phi^\mathrm{E}}}\mathrm{D}_\mathrm{T}\mathcal{J}_{pqk}\left(\mathbf{x}\right),
\label{singlefrequency_td3d}
\end{eqnarray}
and the multi-frequency topological derivative  is defined as the weighted sum:
\begin{eqnarray}
\mathrm{D}_\mathrm{T}\mathcal{J}\left(\mathbf{x}\right) = \frac{1}{N_\mathrm{freq}}\sum_{k=1}^{N_\mathrm{freq}}\frac{\mathrm{D}_\mathrm{T}\mathcal{J}_{k}\left(\mathbf{x}\right)}{\vert\min_{\mathbf{y}\in{\cal R}^{\rm insp}}\mathrm{D}_{\mathrm{T}}\mathcal{J}_{k}\left(\mathbf{y}\right)\vert}.
\label{multifrequency_td3d}
\end{eqnarray}
Also, as in the 2D case, targets will be approximated by considering the sets (directed on a tunnable parameter $\lambda$) of all points in the inspection region for which the topological derivative attains the largest negative values  defined in (\ref{eq:dom_app}).

We recall that the topological energy can be computed  for each experiment as
\begin{equation}\label{eq:ETmaxwell}
\mathrm{E}_\mathrm{T}\mathcal{J}_{pqk}\left(\mathbf{x}\right) = \vert\boldsymbol{\mathcal{U}}_{pqk}\left(\mathbf{x}\right)\vert^2\vert\boldsymbol{\mathcal{V}}_{pqk}\left(\mathbf{x}\right)\vert^2.
\end{equation}
Combinations of experiments and frequencies are done analogously to the topological derivative combinations  and targets are approximated by (\ref{eq:omegaap2}).

\subsection{Target reconstruction under parallel polarization.}
\label{sec:target3d}

\begin{figure}[hbt!]
\centering
\includegraphics[width=10cm]{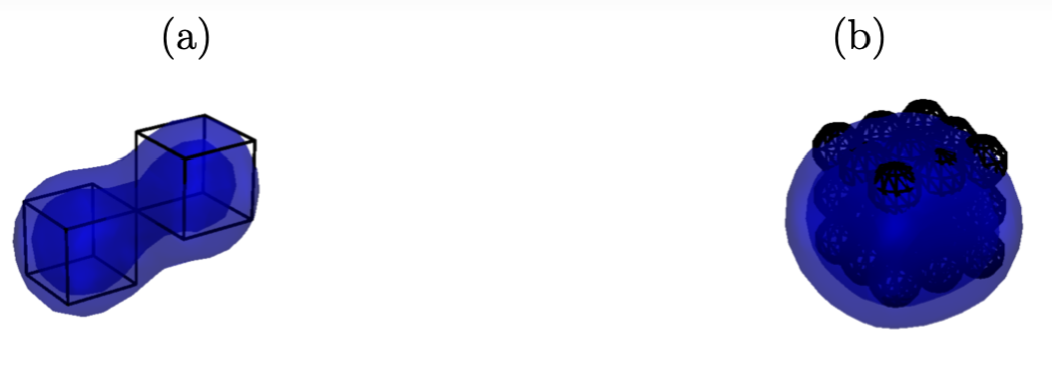}
\caption{Approximations  (\ref{eq:dom_app}) obtained with $\lambda=0.5$ (light) and $\lambda=0.7$ (dark) for (a) the two cubes and (b)  the cube of spheres when the multi-frequency topological derivative (\ref{multifrequency_td3d}) is used on the PP datasets. Black lines represent the true target.}
\label{fig:two_cubes_DTPP_multi}
\end{figure}

Figure \ref{fig:two_cubes_DTPP_multi} represents approximations  $\Omega_{\rm app,\lambda}$ defined in (\ref{eq:dom_app}) obtained for $\lambda=0.5$ and $0.7$ using the parallel polarization dataset of the targets consisting in two cubes and a cube of spheres. In panel (a), we see a certain approximation of the shape of the cubes, but we cannot distinguish if there is one or two objects in the scene. The size and position, nevertheless, are well approximated for a one-step method which uses no a priori information. Similar conclusions can be deduced from panel (b). The topological derivative accurately shows the position and approximate size of the object. However, it is unable to distinguish every element in the aggregate.

\begin{figure}[hbt!]
\centering
\includegraphics[width=14cm]{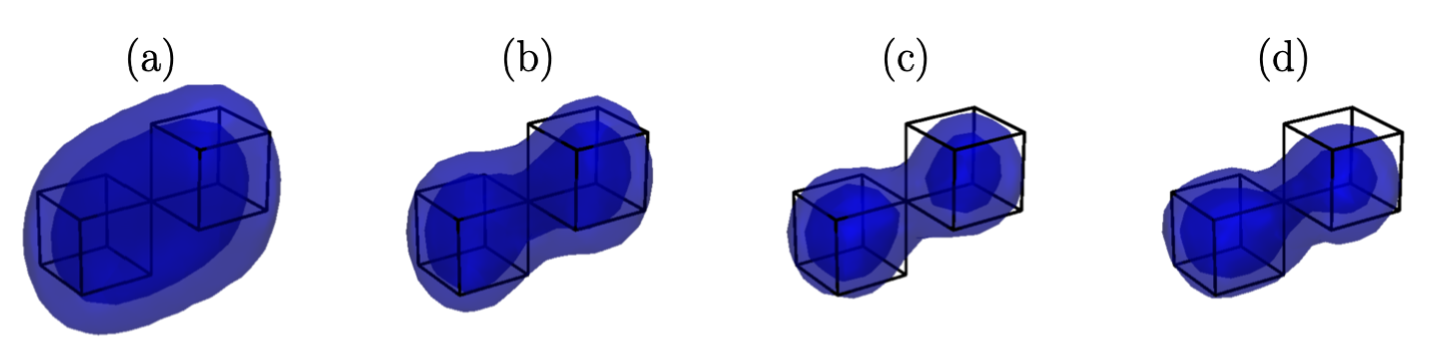}
\caption{Approximations (\ref{eq:dom_app})   for $\lambda=0.5$ (light) and $\lambda=0.7$ (dark) using the single-frequency topological derivative (\ref{singlefrequency_td3d}) with data collected for the frequencies: (a) 3, (b) 4.25, (c) 5.5 and (d) 6.75 GHz. Black lines represent the true target.}
\label{fig:two_cubes_DTPP_freqs}
\end{figure}

\begin{figure}[hbt!]
\centering
\includegraphics[width=14cm]{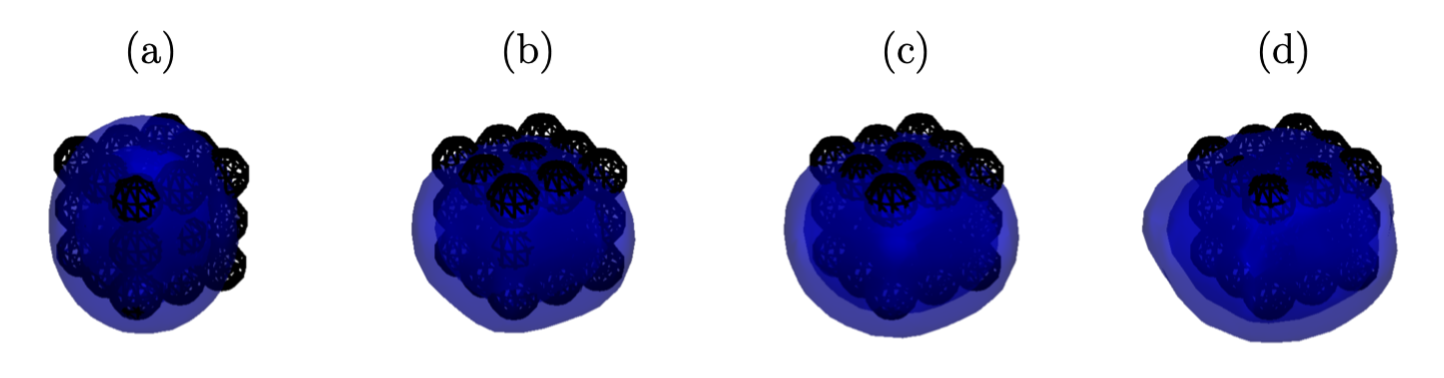}
\caption{ Counterpart of Figure \ref{fig:two_cubes_DTPP_freqs} for the cube of spheres. The considered frequencies are: (a) $3$, (b) $4.25$, (c) $5.5$ and (d) $6.75$  GHz. }\label{fig:cube_spheres_DTPP_freqs}
\end{figure}

The following examples correspond to the smallest targets in the dataset. For these targets, all the frequencies yield topological derivatives with useful information. In fact, some
isolated frequencies provide better information on the target than the multifrequency approach, see panel \ref{fig:two_cubes_DTPP_freqs} (c) in which two separate objects
are identified. This figure illustrates approximations of the cubes provided by single-frequency topological derivatives for different frequencies.  At lower frequencies the topological derivative is less sensitive to the shape variations. A similar phenomenon is observed in Figure \ref{fig:cube_spheres_DTPP_freqs} where the cube of spheres is considered.

\begin{figure}[hbt!]
\centering
\includegraphics[width=14cm]{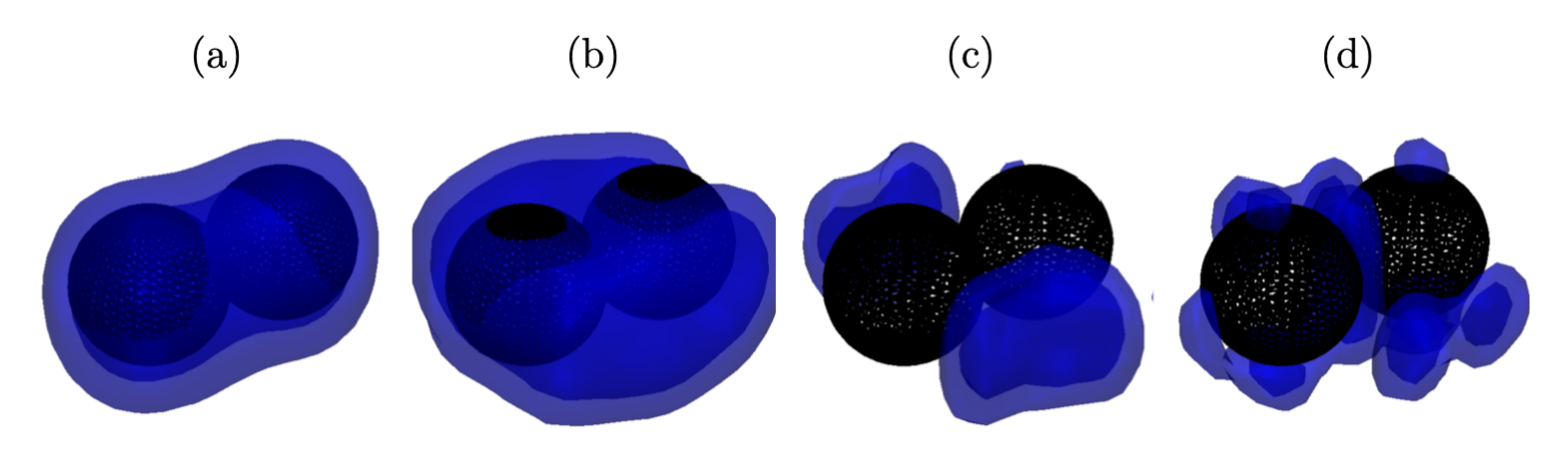}
\caption{ Counterpart of Figure \ref{fig:two_cubes_DTPP_freqs} for the two spheres. The considered frequencies are: (a) $3$, (b) $4.25$, (c) $5.5$ and (d) $6.75$  GHz.}
\label{fig:two_spheres_DTPP_freqs}
\end{figure}

\begin{figure}[hbt!]
\centering
\includegraphics[width=14cm]{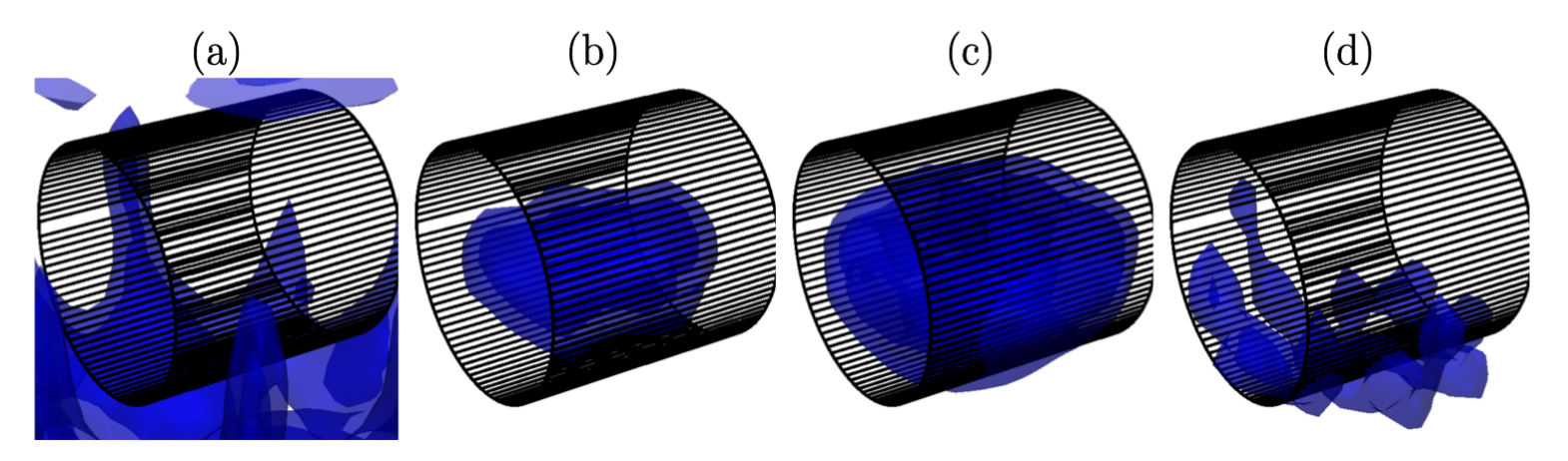}
\caption{Counterpart of Figure \ref{fig:two_cubes_DTPP_freqs} for the cylinder. The considered frequencies in this case are: (a) $4.25$, (b) $5.5$ (c) $6.75$ and (d) $8$ GHz.
}
\label{fig:cylinder_DTPP_freqs}
\end{figure}

For the two largest targets, the cylinder and the two spheres  we did not get satisfactory multi-frequency results.  Being much bigger than the previous target components, only  topological derivatives for some frequencies provide relevant information about them. In Figure \ref{fig:two_spheres_DTPP_freqs} we represent single-frequency topological derivatives applied to the two spheres dataset. While for low frequencies the approximations are reasonably good, as we increase the frequency their quality worsens.
For the cylinder reasonable approximations are only obtained for intermediate values of the frequency, see Figure \ref{fig:cylinder_DTPP_freqs}.

Being a one-step method based on the evaluation of some closed-form expressions, it is much faster  than any of the iterative methods  that were tested in the special sesion. Using a $4$ cores laptop, each target required $3$ to $4$ minutes to perform the multi-frequency reconstructions.

No results are shown for the transverse polarization cases, as we were not able to extract successful reconstructions in a reliable manner for any of the targets. This fact is consistent with observations made in \cite{chaumet2009three} or \cite{geffrin2009continuing}. The targets not having a strong depolarizing effect implies that the transverse polarization data  has very small values, sometimes in the magnitude of the measurement error.

\subsection{Topological energy and reciprocity theorem.}
\label{sec:te3d}

\begin{figure}[h!]
\centering
\includegraphics[width=10cm]{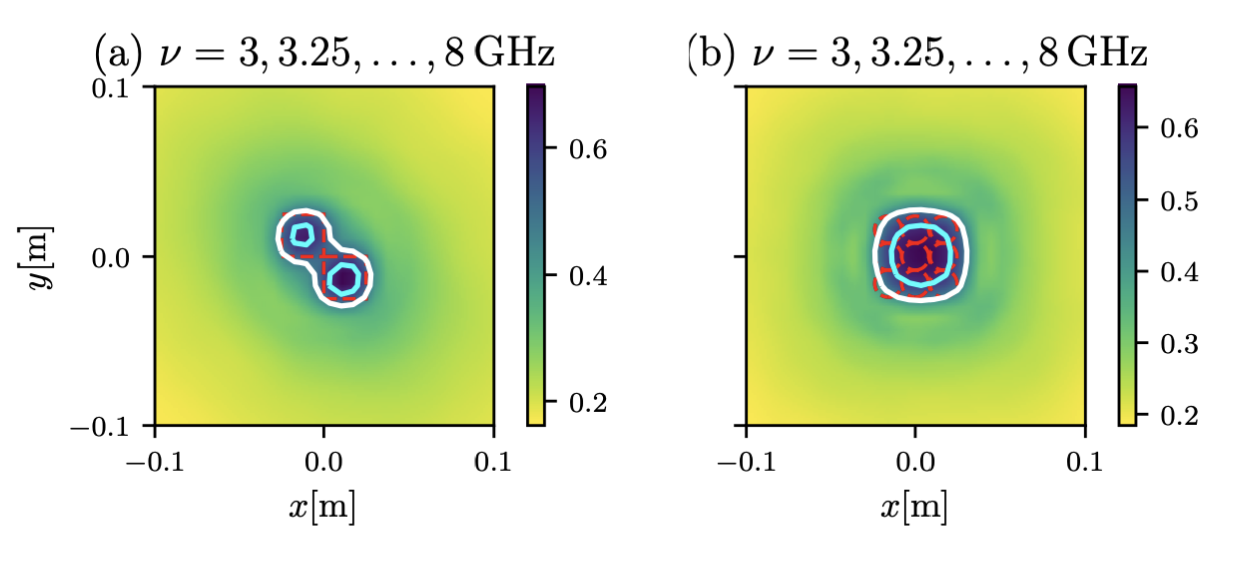}
\caption{Multi-frequency topological energy on the horizontal plane with parallel polarization for (a) the two cubes and (b) the cube of spheres, the smallest targets.
Dashed red lines mark the true targets while solid white  and cyan contours
represent the levels $\lambda=0.7$ and $\lambda=0.9$, respectively.
}
\label{fig:E_twocubes_cubes_spheres}
\end{figure}

\begin{figure}[h!]
\centering
\includegraphics[width=10cm]{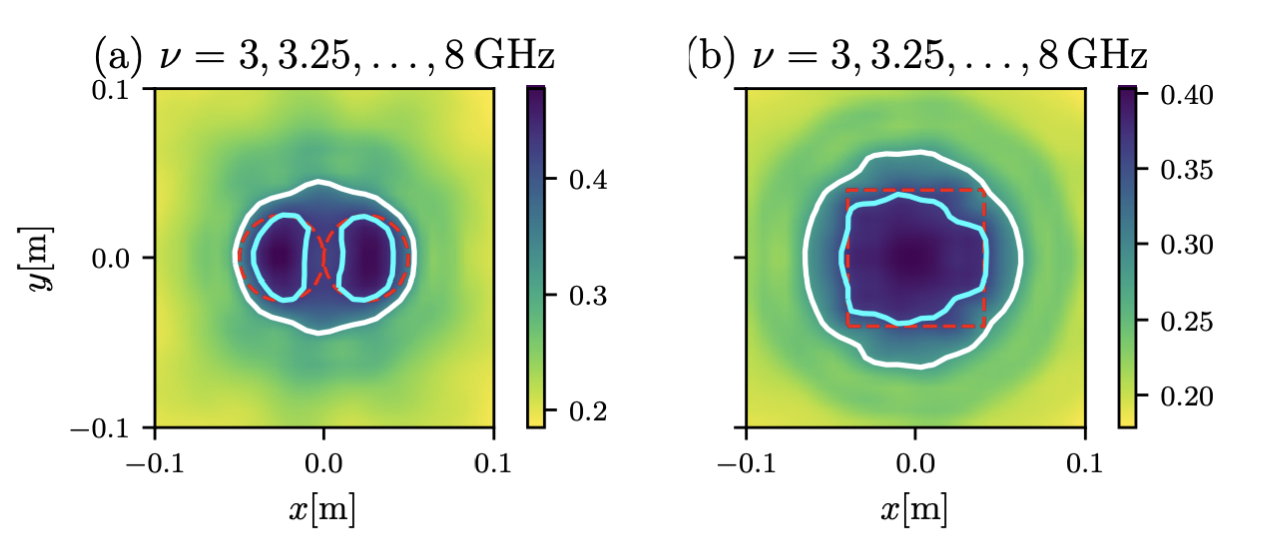}
\caption{Multi-frequency topological energy on the horizontal plane with parallel polarization for (a) the two spheres and (b) the cylinder, the largest targets.
Dashed red lines mark the true targets while solid white  and cyan contours
represent the levels $\lambda=0.7$ and $\lambda=0.9$, respectively.}
\label{fig:E_twospheres_cylinder}
\end{figure}

In three dimensions the topological energy produces  completely different results than those obtained for the topological derivative method.
If we compute the multi-frequency topological energy for  parallel polarization we get a field which is independent of $z$. This is due to the fact that measurements are
performed in the horizontal plane along the vertical direction, and the topological energy formula  (\ref{eq:ETmaxwell}) does not include information on the phase but only on the amplitude.
Figures \ref{fig:E_twocubes_cubes_spheres} and \ref{fig:E_twospheres_cylinder}
represent the multi-frequency topological energy on the horizontal plane for the four  targets in Figure \ref{fig:3d_targets}. That slice of the topological energy is a very accurate indicator of the shape and position of the object, not having spurious minima.

As mentioned in \cite{geffrin2009continuing},  we can use the reciprocity theorem. The reciprocity theorem applied to Maxwell's equations says that, if we have an experiment $\mathrm{A}$ where a current distribution $\mathbf{J}_\mathrm{A}$ causes an electric field $\mathbf{E}_\mathrm{A}$ and another experiment $\mathrm{B}$ where a current distribution $\mathbf{J}_\mathrm{B}$ causes an electric field $\mathbf{E}_\mathrm{B}$, then following equality holds:
 \[
\int_{ \mathbb R^3} \mathbf{J}_\mathrm{A}\cdot\mathbf{E}_\mathrm{B}\,\mathrm{d}V=\int_{\mathbb R^3} \mathbf{J}_\mathrm{B}\cdot\mathbf{E}_\mathrm{A}\,\mathrm{d}V.
\]
The above identity particularized for these experiments, considering the current a point dipole, reads:
\[
\boldsymbol{\mathcal{E}}\left({\bf x}^\mathrm{R};\mathbf{x}^\mathrm{E},\mathbf{p}^\mathrm{E}\right)\cdot\mathbf{k} = \boldsymbol{\mathcal{E}}\left(\mathbf{x}^\mathrm{E};\mathbf{x}^\mathrm{R},\mathbf{k}\right)\cdot\mathbf{p}^\mathrm{E},
\]
where $\boldsymbol{\mathcal{E}}\left(\mathbf{x};\mathbf{x}_0,\mathbf{p}\right)$ is the electric field at the point $\mathbf{x}$ produced by a current dipole $\mathbf{p}$ situated at $\mathbf{x}_0$.
For  the 3D Fresnel database, this implies that we can change the role of emitter and receiver antennas, that is, we can consider any measurement as the one that will be obtained when the object is irradiated vertically polarized from an antenna situated at the measuring point, and this field is measured at point $\mathbf{x}^\mathrm{E}$ with an antenna polarized along $\mathbf{u}_\phi$.  For the transverse polarization case it would be $\mathbf{u}_\theta$ instead of $\mathbf{u}_\phi$, but as already mentioned, this case is not considered due to the amount of noise in the measurements.


\begin{figure}[hbt!]
\centering
\includegraphics[width=14cm]{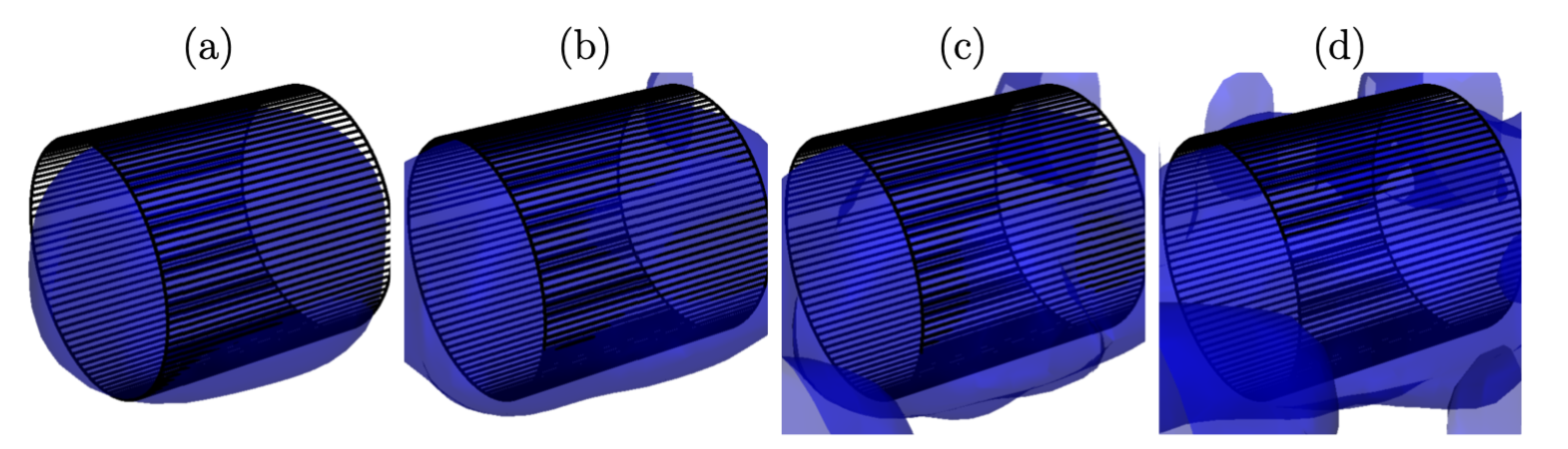}
\caption{Approximations for $\lambda=0.6$  using the single-frequency topological energy with data collected  for the frequencies: (a) $3$, (b) $3.25$ (c) $3.5$ and (d) $3.75$ GHz. Black lines represent the true target.}
\label{fig:cylinder_ETPP_rec_freqs}
\end{figure}

\begin{figure}[hbt!]
\centering
\includegraphics[width=14cm]{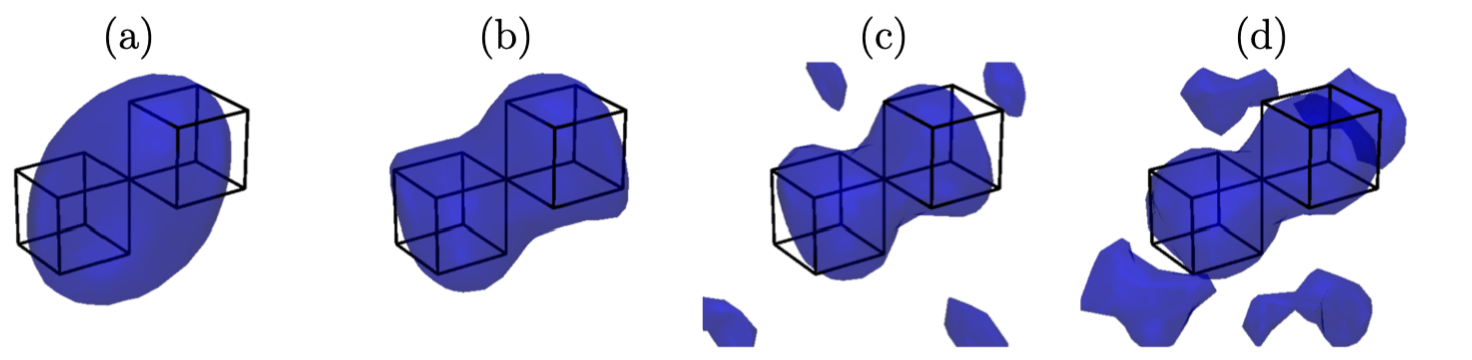}
\caption{Counterpart of  Figure \ref{fig:cylinder_ETPP_rec_freqs} for the two cubes.
 The considered frequencies are:  (a) $3$, (b) $4.25$, (c) $5.5$, (d) $6.75$.}\label{fig:two_cubes_ETPP_rec_freqs}
\end{figure}

\begin{figure}[hbt!]
\centering
\includegraphics[width=14cm]{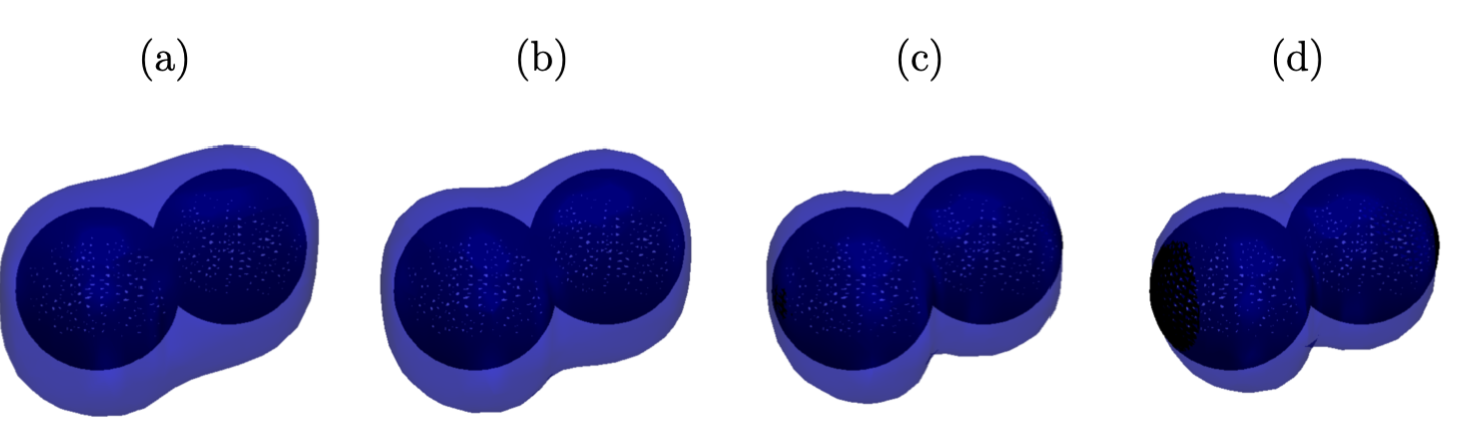}
\caption{ Counterpart of  Figure \ref{fig:cylinder_ETPP_rec_freqs} for the two spheres. The considered frequencies are:  (a) $3$, (b) $3.25$, (c) $3.5$, (d) $3.75$. }
\label{fig:two_spheres_ETPP_rec_freqs}
\end{figure}

\begin{figure}[hbt!]
\centering
\includegraphics[width=14cm]{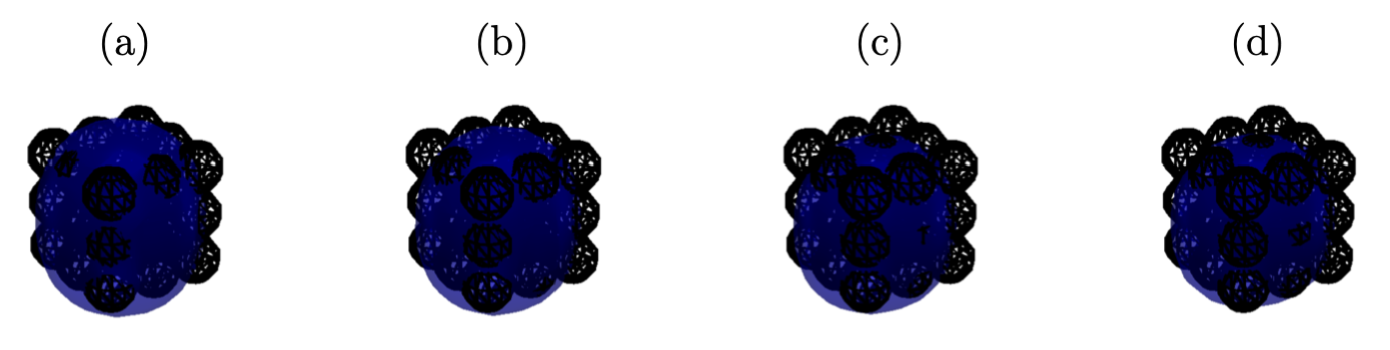}
\caption{Counterpart of  Figure \ref{fig:cylinder_ETPP_rec_freqs} for the cube of spheres.   The considered frequencies are: (a) $3$, (b) $3.25$, (c) $3.5$, (d) $3.75$. 
}
\label{fig:cube_spheres_ETPP_rec_freqs}
\end{figure}


The results obtained with the topological  energy  applied to the data with the emitting and receiving roles interchanged are no longer $z$-independent. However they are very similar to   their counterpart results  obtained with the topological derivative. In this case, the two cubes are the only target for which the approximations are somehow independent of the frequency. For the remaining targets the high frequencies give  highly unsatisfactory  results. Figures \ref{fig:cylinder_ETPP_rec_freqs}, \ref{fig:two_cubes_ETPP_rec_freqs}, \ref{fig:two_spheres_ETPP_rec_freqs} and \ref{fig:cube_spheres_ETPP_rec_freqs} display results for  single-frequency topological energies. Only the frequencies in the range that yields appropriate reconstructions are shown for each target. We select the  value  $\lambda=0.6$ (intermediate between the two used with  the topological derivative), for the ease of visualization.
\section{Conclusions}
\label{sec:conclusions}

Topological imaging methods applied to the 2D Fresnel database provide
a robust and efficient technique yielding good quality approximations to the
targets, with an accuracy comparable to the results obtained with the
distorted-wave Born method   \cite{tijhuis2001multiple} or the contrast source
inversion method  \cite{bloemenkamp2001inversion}. As it happens with those
methods, the topological fields do not need a priori information and are able
to successfully combine multi-frequency information. Being able to combine
all the frequencies means that we can expect  accurate  reconstructions
without knowing in advance which frequencies are best suited for each target.

Being a one step method, topological derivative/energy based imaging
is in principle faster than the iterative techniques employed in
\cite{ baussard2001bayesian,belkebir2001modified2, bloemenkamp2001inversion, crocco2001inverse, duchene2001inversion, fatone2001image,  marklein2001linear, ramananjaona2001shape}.
It  is not possible to perfom a fair comparison of the computational time, as
there are about $20$ years of difference in the hardware used. With a standard
$4$ core desktop computer we can get multi-frequency reconstructions in about
$10$ seconds per target. However, many of these iterative methods address a
more demanding problem,  since they aim not only to recover the shape of the
targets, but also the contrast value (their electrical permittivity or conductivity).
If the actual value of the contrast parameters was needed we could combine the topological derivative method with some iterative scheme for the additional
parameter, as was done in \cite{carpio2012hybrid} with gradient methods with
good results.
Iterative topological derivative based methods \cite{carpio2012hybrid} or hybrid strategies combining Newton approaches for optimizing shapes \cite{carpio2019holography} or parameters \cite{carpio2020bayesian} with
topological techniques can improve the results further.

If we compare with the non-iterative scheme proposed in \cite{testorf2001imaging},
we find that although this method is able to recover somehow the contrast
parameters, they need a priori information on the location of the targets and the
approximated objects obtained are not so close to the actual shapes.

When only the shape of the scatterers is to be found from this database, the
topological derivative method alone is an efficient technique, as it recovers the
shapes of the objects  in a fast way, and without the need of a priori information
of any kind. Moreover, the topological energy provides sharper  maxima (that is,
more robust reconstructions),  requires less information  and if the incident wave
can be modelled as a plane wave it does only need the measured values of the
 scattered field at the receiving antennas. In agreement with previous work,
we observe small deviations in the position or angle in certain configurations,
which may reflect the effect of noise in
the data or aberrations in the imaging system. More careful analyses in the line of
\cite{carpio2020bayesian} could shed light on the role of noise, aberrations and
errors introduced by the approximation of the 3D Maxwell equations by a 2D
Maxwell problem \cite{carpio2018experimental}.

The 3D Fresnel database is more challenging, and the choice of frequency ranges
becomes a key issue. The topological energy and topological derivative methods provide reasonable results on some of the datasets of the 3D database, combining
all frequencies or for specific frequencies. An automatic method for choosing the appropriate frequency range should be developed to proceed in a completely
automatic way, as it was the case for the 2D database.

When processing the parallel polarized datasets, the topological energy gives $z$-independent solutions, similar to the linear sampling method  \cite{catapano20093d} and needs the reciprocity theorem  in order to give full 3D approximations. This does not happen with the topological derivative, which is able to approximate the 3D shape
of the targets using only the parallel polarized data. In this sense, it seems that the topological derivative is more capable of  exploiting the physical symmetries of the underlying model.

Compared to tests performed with other methods on this database  \cite{catapano20093d,chaumet2009three, de2009three, eyraud2009microwave, li2009application, yu2009reconstruction}, the quality of the results is  similar. The
problem of choosing the  suitable frequency is common to all the methods.
Since it is a one-step method that only requires very inexpensive computations,
it is faster than the iterative ones. Using a $4$ cores laptop,  it takes  $3$ to $4$ minutes to perform the multi-frequency reconstruction of each target. It  is not
possible to do a fair comparison of the computation time  as
there are about $10$ years of difference in the hardware used.

As future work, finding a reliable way to choose the threshold $\lambda$ would allow
to generate approximations in an automatic way rather than trying different values of
$\lambda$. It would also be interesting to study more advanced ways of combining
the different  single-frequency topological derivatives and selecting the most
informative frequency ranges. Exploring the performance of iterative methods based
on the computation of iterated topological derivatives in the spirit of \cite{carpio2008topological} could also provide better reconstructions
but at a much higher computational cost.

\section*{Acknowledgements}
Research partially supported by Spanish FEDER/MICINN-AEI grants MTM2014-56948-C2-1-P, MTM2017-84446-C2-1-R and TRA2016-75075-R.


\end{document}